\documentclass{amsart}
\usepackage{amsmath,amsthm}
\usepackage{amssymb,amscd,enumerate}
\usepackage{float,epsfig}
\usepackage{pb-diagram}
\def\a{\alpha}
\def\b{\beta}

\def\d{\delta}

\def\e{\epsilon}

\def\G{\Gamma}

\def\k{\kappa}
\def\lb\{{\left\{}
\def\la{\lambda}
\def\La{\Lambda}
\def\lla{\longleftarrow}
\def\lm{\limits}
\def\lra{\longrightarrow}
\def\dllra{\Longleftrightarrow}
\def\llra{\longleftrightarrow}
\def\n{\nabla}
\def\ngth{\negthickspace}
\def\ola{\overleftarrow}
\def\Om{\Omega}
\def\om{\omega}
\def\op{\oplus}
\def\oper{\operatorname}
\def\oplm{\operatornamewithlimits}
\def\ora{\overrightarrow}
\def\ov{\overline}
\def\ova{\overarrow}
\def\ox{\otimes}
\def\p{\partial}
\def\rb\}{\right\}}
\def\s{\sigma}
\def\sbq{\subseteq}
\def\spq{\supseteq}
\def\sqp{\sqsupset}
\def\supth{{\text{th}}}
\def\T{\Theta}
\def\th{\theta}
\def\tl{\tilde}
\def\thra{\twoheadrightarrow}
\def\un{\underline}
\def\ups{\upsilon}
\def\vp{\varphi}
\def\wh{\widehat}
\def\wt{\widetilde}
\def\x{\times}
\def\z{\zeta}
\def\({\left(}
\def\){\right)}
\def\[{\left[}
\def\]{\right]}
\def\<{\left<}
\def\>{\right>}

\def\tec{Teichm\"uller\ }
\def\sconr{\hbox{\medspace\vrule width 0.4pt height 4.7pt depth
0.4pt \vrule width 5pt height 0pt depth 0.4pt\medspace}}

\def\SA{\mathcal A}
\def\SB{\mathcal B}
\def\SC{\mathcal C}
\def\SD{\mathcal D}
\def\SE{\mathcal E}
\def\SF{\mathcal F}
\def\SG{\mathcal G}
\def\SH{\mathcal H}
\def\SI{\mathcal I}
\def\SJ{\mathcal J}
\def\SK{\mathcal K}
\def\SL{\mathcal L}
\def\SM{\mathcal M}
\def\SN{\mathcal N}
\def\SO{\mathcal O}
\def\SP{\mathcal P}
\def\SQ{\mathcal Q}
\def\SR{\mathcal R}
\def\SS{\mathcal S}
\def\ST{\mathcal T}
\def\SU{\mathcal U}
\def\SV{\mathcal V}
\def\SW{\mathcal W}
\def\SX{\mathcal X}
\def\SY{\mathcal Y}
\def\SZ{\mathcal Z}

\newcommand{\BA}{\ensuremath{\mathbf A}}
\newcommand{\BB}{\ensuremath{\mathbf B}}
\newcommand{\BC}{\ensuremath{\mathbf C}}
\newcommand{\BD}{\ensuremath{\mathbf D}}
\newcommand{\BE}{\ensuremath{\mathbf E}}
\newcommand{\BF}{\ensuremath{\mathbf F}}
\newcommand{\BG}{\ensuremath{\mathbf G}}
\newcommand{\BH}{\ensuremath{\mathbf H}}
\newcommand{\BI}{\ensuremath{\mathbf I}}
\newcommand{\BJ}{\ensuremath{\mathbf J}}
\newcommand{\BK}{\ensuremath{\mathbf K}}
\newcommand{\BL}{\ensuremath{\mathbf L}}
\newcommand{\BM}{\ensuremath{\mathbf M}}
\newcommand{\BN}{\ensuremath{\mathbf N}}
\newcommand{\BO}{\ensuremath{\mathbf O}}
\newcommand{\BP}{\ensuremath{\mathbf P}}
\newcommand{\BQ}{\ensuremath{\mathbf Q}}
\newcommand{\BR}{\ensuremath{\mathbf R}}
\newcommand{\BS}{\ensuremath{\mathbf S}}
\newcommand{\BT}{\ensuremath{\mathbf T}}
\newcommand{\BU}{\ensuremath{\mathbf U}}
\newcommand{\BV}{\ensuremath{\mathbf V}}
\newcommand{\BW}{\ensuremath{\mathbf W}}
\newcommand{\BX}{\ensuremath{\mathbf X}}
\newcommand{\BY}{\ensuremath{\mathbf Y}}
\newcommand{\BZ}{\ensuremath{\mathbf Z}}

\def\bba{{\mathbb A}}
\def\bbb{{\mathbb B}}
\def\bbc{{\mathbb C}}
\def\bbd{{\mathbb D}}
\def\bbe{{\mathbb E}}
\def\bbf{{\mathbb F}}
\def\bbg{{\mathbb G}}
\def\bbh{{\mathbb H}}
\def\bbi{{\mathbb I}}
\def\bbj{{\mathbb J}}
\def\bbk{{\mathbb K}}
\def\bbl{{\mathbb L}}
\def\bbm{{\mathbb M}}
\def\bbn{{\mathbb N}}
\def\bbo{{\mathbb O}}
\def\bbp{{\mathbb P}}
\def\bbq{{\mathbb Q}}
\def\bbr{{\mathbb R}}
\def\bbs{{\mathbb S}}
\def\bbt{{\mathbb T}}
\def\bbu{{\mathbb U}}
\def\bbv{{\mathbb V}}
\def\bbw{{\mathbb W}}
\def\bbx{{\mathbb X}}
\def\bby{{\mathbb Y}}
\def\bbz{{\mathbb Z}}

\newtheorem{thm}{Theorem}[section]
\newtheorem{lem}[thm]{Lemma}
\newtheorem{prop}[thm]{Proposition}
\newtheorem{cor}[thm]{Corollary}

\theoremstyle{definition}
\newtheorem{defn}[thm]{Definition}
\newtheorem{rem}[thm]{Remark}

\newcommand{\1}{\ensuremath{\pi_1}}
\newcommand{\2}[1]{\ensuremath{^{(#1)}}}
\newcommand{\nn}{\ensuremath{^{(n+1)}}}
\newcommand{\ff}{\ensuremath{\SF_{(n.0)}/\SF_{(n.5)}}}
\newcommand{\zz}{\ensuremath{\SG_{(n+2)}/\SG_{(n+2.5)}}}
\newcommand{\sr}{\ensuremath{S^3\backslash R}}
\newcommand{\sk}{\ensuremath{S^3\backslash K}}
\newcommand{\ns}{$n$-solvable}
\newcommand{\nss}{$(n.5)$-solvable}

\newcommand{\arf}{\operatorname{Arf}}
\newcommand{\Hom}{\operatorname{Hom}}
\newcommand{\id}{\operatorname{id}}
\newcommand{\image}{\operatorname{image}}
\newcommand{\Image}{\operatorname{Image}}
\newcommand{\kernel}{\operatorname{kernel}}
\newcommand{\Ker}{\operatorname{Ker}}
\newcommand{\rank}{\operatorname{rank}}
\newcommand{\inc}{\operatorname{inc}}
\newcommand{\rk}{\operatorname{rk}}
\newcommand{\PD}{\operatorname{PD}}

\def\kg{{\SK_\G}}
\def\ee{\acute{\text{e}}}
\def\Anz{{\SA^\bbz_n}}
\def\An{{\SA_n}}
\def\tg{\tilde{G}}
\def\br{\b_\text{rel}}

\begin{document}

\title{Higher-order Alexander invariants and
filtrations of the knot concordance group}
\author{Tim D. Cochran and Taehee Kim $^\dag$}
\date{\today}
\address{Rice University, Houston, Texas, 77005-1892}
\email{cochran@math.rice.edu, tkim@rice.edu}
\thanks{$^\dag$ The first author was partially supported
by the National Science Foundation}

\begin{abstract} We establish certain ``nontriviality'' results
for several filtrations of the smooth and topological knot
concordance groups. First, as regards the \emph{n-solvable
filtration} of the topological knot concordance group, $\SC$,
defined by K. Orr, P. Teichner and the first author \cite{COT1}:
$$
 0\subset\cdots\subset\SF_{(n.5)}\subset\SF_{(n)}\subset\cdots
 \subset\SF_{(1.5)}\subset\SF_{(1.0)}\subset\SF_{(0.5)}
 \subset\SF_{(0)}\subset\SC,
$$
we refine the recent nontriviality results of Cochran and Teichner
\cite{CT} by including information on the Alexander modules. These
results also extend those of C. Livingston \cite{Li} and the
second author \cite{K}. We exhibit similar structure in the
closely related \emph{symmetric Grope filtration} of $\SC$
considered in \cite{CT}. We also show that the Grope filtration of
the \emph{smooth} concordance group is nontrivial using examples
that cannot be distinguished by the Ozsv{\'a}th-Szab{\'o}
$\tau$-invariant nor by J. Rasmussen's $s$-invariant
\cite{OS}\cite{Ra}. Our broader contribution is to establish, in
``the relative case'', the key homological results whose analogues
Cochran-Orr-Teichner established in ``the absolute case'' in
\cite{COT1}.

We say two knots $K_0$ and $K_1$ are \emph{concordant modulo
$n$-solvability} if $K_0\#(-K_1)\in \SF_{(n)}$. Our main result is
that, for any knot $K$ whose classical Alexander polynomial has
degree greater than 2, and for any positive integer $n$, there
exist infinitely many knots $K_i$ that are concordant to $K$
modulo $n$-solvability, but are all distinct modulo
$n.5$-solvability. Moreover, the $K_i$ and $K$ share the same
classical Seifert matrix and Alexander module as well as sharing
the same \emph{higher-order} \emph{Alexander modules} and Seifert
presentations up to order $n-1$.

\end{abstract}

\baselineskip=13pt \maketitle
\section{Introduction}
Oriented knots $K_0$ and $K_1$ in the 3-sphere are called
(topologically) \emph{concordant} if there is a topologically
locally flat proper embedding of $S^1\times [0,1]$ into $S^3\times
[0,1]$ that restricts to the $K_i$ on $S^3\times \{i\}$.
Equivalently, $K_0$ is concordant to $K_1$ if $K_0\#(-K_1)$ is
(topologically) slice, that is, if $K_0\#(-K_1)$ bounds a
topologically locally flat disk in the 4-ball. Here $`\#'$ denotes
the connected sum, and $-K_1$ is the mirror image of $K_1$ with
reversed string orientation. Concordance is an equivalence
relation on the set of oriented knots, and the set of equivalence
classes (\emph{concordance classes}) forms an abelian group,
$\SC$, under the operation of connected sum. This group is called
the \emph{(topological) knot concordance group}. The
classification of the knot concordance group is still open, and it
has been one of the central problems in knot theory. A common
strategy of investigating the knot concordance group is to extract
information from abelian covers or metabelian covers of the
exterior of a knot (see \cite{L,CS,CG,G,KL,Le,Fr}). On the other
hand, requiring smooth embedding instead of topologically locally
flat embedding, we can define the \emph{smooth knot concordance
group}, which is denoted by $\SC^{smooth}$. In smooth knot
concordance theory, there have been recent and very interesting
developments using knot Floer homology \cite{OS} and Khovanov
homology \cite{Ra}. In this paper we concern ourselves mainly with
topological knot concordance. Hence we work in the topologically
locally flat category unless mentioned otherwise.

Cochran, Orr, and Teichner (henceforth COT) recently made
significant progress in the study of $\SC$ using \emph{derived}
covers of the exterior of a knot \cite{COT1}. (A \emph{derived
cover} is a covering space corresponding to a derived subgroup of
the fundamental group of the space.) In particular, they defined a
filtration $\{\SF_{(n)}\}_{n\in \frac12 \bbn_0}$ of the knot
concordance group:
$$
 0\subset\cdots\subset\SF_{(n.5)}\subset\SF_{(n)}\subset\cdots
 \subset\SF_{(1.5)}\subset\SF_{(1.0)}\subset\SF_{(0.5)}
 \subset\SF_{(0)}\subset\SC.
$$
Here $\SF_{(n)}$ is a subgroup of $\SC$ consisting of
\emph{$n$-solvable} knots. (A knot $K$ is said to be
\emph{$n$-solvable} if the zero surgery on the knot in the
3-sphere bounds a spin 4-manifold $W$ which satisfies certain
conditions on integral homology groups and intersection form on
the $n$-th derived cover of $W$. In this case, we say $K$ is
\emph{$n$-solvable via $W$}, and $W$ is called an
\emph{$n$-solution} for $K$. Refer to \cite[Section 8]{COT1}.) COT
showed that the previously known abelian and metabelian
concordance invariants are reflected at the first stage of their
filtration \cite{COT1}, implying that the abelian group
$\SF_{(1.0)}/\SF_{(1.5)}$ has infinite rank. They also established
that $\SF_{(2.0)}/\SF_{(2.5)}$ has infinite rank \cite{COT2}.
Moreover, Cochran and Teichner showed that $\SF_{(n)}/\SF_{(n.5)}$
has positive rank for every integer $n \ge 2$ \cite{CT}. It is
still unknown whether or not $\SF_{(n)}/\SF_{(n.5)}$ is infinitely
generated for $n>2$, or whether or not $\SF_{(n.5)}/\SF_{(n+1)}$
is non-trivial. Thus, although $\SF$ has been shown to be highly
non-trivial, many questions about its structure remain open.

These non-triviality results have been refined in the case $n=1$.
C. Livingston recently asked if one can always find non-concordant
knots which share a given (classical) Seifert form. That is, given
$K$, can one always find other knots, distinct up to concordance,
which share the same (classical) Seifert form as $K$? As he
observed, if $K$ has the Seifert form of an unknot (has Alexander
polynomial $1$) then the answer is certainly ``No'', since all
knots with Alexander polynomial $1$ are topologically concordant
to the trivial knot and hence concordant to each other (by M.
Freedman's work \cite{F,FQ}). Livingston gave a partial answer in
the positive to his question using Casson-Gordon invariants
\cite{Li}. This work was completed by the second author, who
showed that there always exist such knots for a given (classical)
Seifert form if and only if the Alexander polynomial of the
Seifert form is not trivial. More precisely, in \cite{K} he showed
that for a given knot $K$ whose Alexander polynomial has degree at
least $2$, there are infinitely many knots $K_i$ (with $K_0 = K$)
such that $K_i-K$ is $1$-solvable, but $K_i - K_j$ ($i\ne j)$ is
not $(1.5)$-solvable, and $K_i$ shares the same (classical)
Seifert form (hence the same (classical) Alexander module) as the
knot $K$. We view these results as non-triviality results for
$\SF$ that are finer (for $n=1$) than those of
\cite{COT1}\cite{COT2}\cite{CT}.  They suggest further questions.
Is the Livingston-Kim result true for $n>1$? Moreover, is there an
even finer result that constructs such examples while fixing not
only the (classical) Seifert form (hence the (classical) Alexander
module) but also the higher-order Seifert forms and higher-order
Alexander modules developed in \cite{C}? This paper answers these
questions in the positive as long as the degree of the Alexander
polynomial is at least $4$ (see below). To do so we found the need
to extend the technology of \cite{COT1} to a ``relative'' setting,
and to greatly generalize the key technical result of Cochran and
Teichner \cite[Theorem 4.3]{CT}. We expect that this extended
technology will be of independent interest, enabling further work
in this area.

In this paper, we introduce the notion of an \emph{$n$-cylinder}
(Definition~\ref{defn:n-cylinder}), generalizing the notion of an
$n$-solution of COT . The difference is that an $n$-solution
allows only one boundary component, whereas an $n$-cylinder can
have multiple boundary components. A basic (trivial) example of an
$n$-cylinder (with two boundary components) is a (spin) homology
cobordism between the zero surgeries on two knots. If two knots
are concordant then one can easily find such an $n$-cylinder
(which in this case is a homology cobordism), by doing surgery on
$S^3\times [0,1]$ along the annulus cobounded by the knots
(Remark~\ref{rem:n-cylinder}(4)). Using this, we define a family
of new equivalence relations on the knot concordance group, that
we call \emph{$n$-solvequivalence} (see
Definition~\ref{defn:n-solvequivalence}) : $K_0$ is
\emph{$n$-solvequivalent} to $K_1$ if the zero-framed surgery on
$K_0$ and the zero-framed surgery on $K_1$ cobound an
$n$-cylinder. From the basic example above, it is clear that
concordant knots are $n$-solvequivalent for all $n$.

On the other hand, the filtration $\SF$ of COT immediately
suggests another family of equivalence relations on $\SC$, given
as follows. We say $K_0$ and $K_1$ are \emph{concordant modulo
$n$-solvability} if $K_0\#(-K_1)$ is $n$-solvable. There is a
close relationship between $n$-solvability and
$n$-solvequivalence. It is not difficult to show that if two knots
are concordant modulo $n$-solvability, then they are
$n$-solvequivalent (Proposition~\ref{prop:old and new}), but we
have not been able to establish the converse. Hence
$n$-solvequivalence is a possibly weaker obstruction to
concordance than $n$-solvability.
Note that $n$-solvequivalence reflects information on the zero
surgeries on \emph{each} $K_i$ separately, but \emph{concordance
modulo $n$-solvability} reflects only information on the zero
surgery on $K_0\#(-K_1)$ (see Proposition~\ref{prop:P and P
perb}).

As an application of $n$-solvequivalence and of our other
extensions of the technology of \cite{COT1}\cite{COT2} \cite{CT},
we obtain the following theorem, generalizing the nontriviality
results cited above. The condition that the degree of the
Alexander polynomial be at least $4$ may seem ad hoc. However,
very recent work of S. Friedl and Teichner \cite{FT} complements
our theorem and we use it to show that our theorem is ``best
possible", in the sense that it is false for certain knots whose
Alexander polynomial has degree $2$.

The reader is referred to \cite[Section 3]{CT} for the definition
of symmetric Gropes. Let $G$ be a group. Then the \emph{$n$-th
derived group of $G$} is inductively defined by $G^{(0)}\equiv G$
and $G\nn \equiv [G^{(n)},G^{(n)}]$.

\newtheorem*{main}{Theorem~\ref{thm:main}}
\begin{main}[{\bf Main Theorem}]
Let $n$ be a positive integer. Let $K$ be a knot whose Alexander
polynomial has degree greater than 2 (if $n=1$ then degree equal
to $2$ is allowed). Then there is an infinite family of knots
$\{K_i \,\,|\,\, i=0,1,2,\ldots\}$ with $K_0 = K$ that satisfies
the following :

\begin{itemize}
\item [(1)] For each $i$, $K_i - K$ is $n$-solvable. In
particular, $K_i$ is $n$-solvequivalent to $K$. Moreover, $K_i$
and $K$ cobound, in $S^3\times\[0,1\]$, a smoothly embedded symmetric
Grope of height $n+2$.

\item [(2)] If $i\ne j$, then $K_i$ is not $(n.5)$-solvequivalent
to $K_j$. In particular, $K_i-K_j$ is not $(n.5)$-solvable, and
$K_i$ and $K_j$ do not cobound, in $S^3\times\[0,1\]$, an embedded
symmetric Grope of height $(n+2.5)$.

\item [(3)] For each $i$, $K_i$ has the same order $m$ integral
higher-order Alexander module as $K$ for $m=0,1,\ldots, n-1$.
Indeed, if $G_i$ and $G$ denote the knot groups of $K_i$ and $K$
respectively, then there is an isomorphism $G_i/(G_i)^{(n+1)}\to
G/G^{(n+1)}$ that preserves the peripheral structures.

\item [(4)] For each $i$, $K_i$ has the same $m^{th}$-higher-order
Seifert presentation as $K$ for $m=0,1,\ldots, n-1$. In
particular, all of the knots admit the same classical Seifert
matrix.

\item [(5)] If $i>j$, $K_i - K_j$ is of infinite order in
$\SF_{(n)}/\SF_{(n.5)}$.

\item [(6)] If $i,j>0$, $s(K_i)=s(K_j)$ and $\tau(K_i)=\tau(K_j)$,
where $s$ is Rasmussen's smooth concordance invariant and $\tau$
is the smooth concordance invariant of P.~Ozsv{\'a}th and
Z.~Szab{\'o} \cite{Ra}\cite{OS}.

\end{itemize}
\end{main}

The general term \emph{higher-order Alexander modules} was first
used in \cite{COT1}. The \emph{higher-order Alexander modules} and
\emph{higher-order Seifert presentations} of a knot that we use
here were introduced and developed by Cochran \cite{C}. These are
defined in terms of the homology of the derived covers of the knot
exterior. Definitions are given in Section~\ref{sec:main}. But in
particular, these generalize (the case $n=0$) the classical
Alexander module and Seifert matrix. Therefore, the case $n=1$ of
our theorem recovers the previously mentioned results of
Livingston and Kim. If the degree of the Alexander polynomial is
precisely $2$ then the Livingston-Kim result says that it is
possible to fix the \emph{classical} Seifert matrix ($m=0$), but
we show, using recent work of Friedl and Teichner, that there are
at least \emph{some} such knots for which it is \emph{not
possible} to fix any higher-order Seifert presentation matrices.
Specifically, the main theorem fails for $n>1$ if we allow knots
whose Alexander polynomials have degree $2$ (see
Proposition~\ref{prop:failure}). More details are within.

As another obvious
corollary we recover the aforementioned result of Cochran and Teichner
\cite{CT}.

\newtheorem*{cor of main}{Corollary~\ref{cor:cot filtration}}
\begin{cor of main}\cite{CT}
For every positive integer $n$, $\SF_{(n)}/\SF_{(n.5)}$ has
positive rank.
\end{cor of main}

\noindent Indeed, this corollary is essentially equivalent to parts (1) and
(5) in Theorem~\ref{thm:main}, and so our new contribution  lies
in being able to impose (3) and (4).

Now consider the filtration  of $\SC$ given by the subgroups
$\SG_n$ of knots that bound topologically embedded symmetric
Gropes of height $n$ in $B^4$ and the filtration of the smooth
knot concordance group $\SC^{smooth}$ given by the subgroups
$\SG_n^{smooth}$ of knots that bound smoothly embedded symmetric
Gropes of height $n$ in $B^4$. It was shown in \cite[Theorem
1.3]{CT} that $\SG_{n+2}/\SG_{n+2.5}$ has positive rank if $n\geq
0$. It follows from their work (although not explicitly stated)
that $\SG_{n+2}^{smooth}/\SG_{n+2.5}^{smooth}$ has positive rank.
Our work recovers these results, together with the added control
of the Alexander and Seifert data, and allows for the following
additional observations concerning the smooth concordance
invariants of Rasmussen and Ozsv{\'a}th-Szab{\'o}.

\newtheorem*{gropecor}{Corollary~\ref{cor:gropefiltrations}}
\begin{gropecor} Let $n$ be a positive integer.
\begin{itemize}
\item [1.]  $\SG_{n+2}^{smooth}/\SG_{n+2.5}^{smooth}$  has
elements of infinite order represented by knots $K$ with $s(K)=0$
and $\tau(K)=0$. \item [2.]  The kernel of the homomorphism
$ST:\SG_{n+2}^{smooth}\to \mathbb{Z}\times \mathbb{Z}$ given by
$ST(K)=(s(K),\tau(K))$ is infinite.
\end{itemize}
\end{gropecor}

Our Proposition~\ref{prop:geomalg} and
Theorem~\ref{thm:paininbutt} represent a significant strengthening
of the crucial technical results of Proposition $6.2$ and Theorem
$6.3$ in \cite{CT}. In the process, we introduce the notion of an
\emph{algebraic $n$-solution} which should be considered as an
algebraic abstraction of our notion of an $n$-cylinder
(Definition~\ref{defn:algebraic} and
Proposition~\ref{prop:geomalg}). In fact, an algebraic
$n$-solution was introduced by Cochran and Teichner in \cite{CT},
but our notion is much more general.

For an obstruction to $n$-solvequivalence, as in
\cite{COT1}\cite{COT2}\cite{CT}, we use the Cheeger-Gromov
\emph{von Neumann $\rho$-invariant} (Theorem~\ref{thm:rho=0} and
Corollary~\ref{thm:rho=0}), \cite{ChG}. For this purpose, we
analyze the twisted homology groups of an $n$-cylinder and its
boundary in Section~\ref{sec:homology}.

This paper is organized as follows. In
Section~\ref{sec:n-cylinder}, we define an $n$-cylinder and
$n$-solvequivalence. We also investigate the relationship between
$n$-solvability and $n$-solvequivalence. In
Section~\ref{sec:homology}, the twisted homology groups of an
$n$-cylinder and its boundary are analyzed. In
Section~\ref{sec:obstruction}, we explain how to realize von
Neumann $\rho$-invariants as obstructions to $n$-solvequivalence.
All of these sections generalize \cite{COT1} by extending to ``the
relative case''. The main theorem is proved in
Section~\ref{sec:main}. In the last section, we define and
investigate an algebraic $n$-solution. Here, we not only
generalize \cite{CT} to the relative case but also significantly
strengthen a primary technical tool. In a strictly logical order,
Section~\ref{sec:algebraic} would precede Section~\ref{sec:main},
but we have chosen to place Section~\ref{sec:algebraic} after
Section~\ref{sec:main}, since the arguments in
Section~\ref{sec:algebraic} are very technical.

\section{$n$-cylinders and $n$-solvequivalence}
\label{sec:n-cylinder} In this section, we define an $n$-cylinder
and $n$-solvequivalence. Throughout this
paper integer coefficients are understood for homology groups
unless specified otherwise.

Recall that a group $G$ is called \emph{\ns} if $G\nn =
1$. If $X$ is a topological space, $X^{(n)}$ denotes the covering
space of $X$ corresponding to the $n$-th derived group of $\1(X)$.
For a compact spin 4-manifold $W$ there is an equivariant
intersection form
$$
\la_n : H_2(W^{(n)}) \times H_2(W^{(n)}) \lra
\bbz[\1(W)/\1(W)^{(n)}]
$$
and a self-intersection form $\mu_n$ \cite[Chapter
5]{Wa}\cite[Section 7]{COT1}. (Here $H_2(W^{(n)})$ is considered
to be a right $\bbz[\1(W)/\1(W)^{(n)}]$-module.) In general, these
intersection forms are singular.

Let us denote the boundary components of $W$ by $M_i$
($i = 1,2,\ldots, \ell$). That is, $\partial W =
\coprod^\ell_{i=1}M_i$ where each $M_i$ is a connected 3-manifold.
Let $I \equiv \Image \{\inc_* : H_2(\partial W) \to H_2(W)\}$
where $\inc_*$ is the homomorphism induced by the inclusion from
$\partial W$ to $W$. Then the usual intersection form factors
through
$$
\ov{\la_0} : H_2(W)/I \times H_2(W)/I \lra \bbz.
$$
We define an \emph{$n$-Lagrangian} to be a
$\bbz[\1(W)/\1(W)^{(n)}]$-submodule of $H_2(W^{(n)})$ on which
$\la_n$ and $\mu_n$ vanish and which maps onto a $\frac12$-rank
direct summand of $H_2(W)/I$ under the covering map. An
\emph{$n$-surface} is defined to be a based and immersed surface
in $W$ that can be lifted to $W^{(n)}$. Observe that any class in
$H_2(W^{(n)})$ can be represented by an $n$-surface and that
$\lambda_n$ can be calculated by counting intersection points in
$W$ among representative $n$-surfaces weighted appropriately by
signs and by elements of $\pi_1(W)/\pi_1(W)^{(n)}$. We say an
$n$-Lagrangian $L$ admits \emph{$m$-duals} (for $m\le n$) if $L$
is generated by (lifts of) $n$-surfaces $\ell_1,\ell_2,\ldots,\ell_g$ and
there exist $m$-surfaces $d_1,d_2,\ldots, d_g$ such that $H_2(W)/I$
has rank $2g$ and $\la_m(\ell_i,d_j)=\delta_{i,j}$. Similarly we
can define a \emph{rational $n$-Lagrangian} and \emph{rational
$m$-duals}. Here we do not require that $W$ be spin and $\mu_n$ is
not discussed. A \emph{rational $n$-Lagrangian} $L$ is a
$\bbq[\1(W)/\1(W)^{(n)}]$-submodule of $H_2(W^{(n)};\bbq)$ on
which $\la_n$ (with $\bbq [\pi_1(W)/\pi_1(W)^{(n)}]$ coefficients)
vanishes and which maps onto a $\frac12$-rank direct summand of
$H_2(W;\bbq)/I_\bbq$ under the covering map. Here $I_\bbq \equiv
\Image \{\inc_* : H_2(\partial W;\bbq) \to H_2(W;\bbq)\}$. We say
a rational $n$-Lagrangian $L$ admits \emph{rational $m$-duals}
(for $m\le n$) if $L$ is generated by $n$-surfaces
$\ell_1,\ell_2,\ldots,\ell_g$ and there are $m$-surfaces
$d_1,d_2,\ldots, d_g$ such that $H_2(W;\bbq)/I_\bbq$ has rank $2g$
and $\la_m(\ell_i,d_j)=\delta_{i,j}$.

\begin{defn}
\label{defn:n-cylinder} Let $n$ be a nonnegative integer.
\begin{enumerate}
\item A compact, connected, spin 4-manifold $W$ with
$\partial W = \coprod^\ell_{i=1}M_i$ where each $M_i$ is a
connected component of $\partial W$ with $H_1(M_i) \cong \bbz$ is
an \emph{$n$-cylinder} if the inclusion from $M_i$ to $W$ induces
an isomorphism on $H_1(M_i)$ for each $i$ and $W$ admits an
$n$-Lagrangian with $n$-duals. In addition to this, if the
$n$-Lagrangian is the image of an $(n+1)$-Lagrangian under the
covering map, then $W$ is called an \emph{$(n.5)$-cylinder}.

\item A compact, connected 4-manifold $W$ with
$H_1(W;\mathbb{Q})\cong \mathbb{Q}$ and with $\partial W =
\coprod^\ell_{i=1}M_i$ where each $M_i$ is a connected component
of $\partial W$ with $H_1(M_i) \cong \bbq$ is a \emph{rational
$n$-cylinder of multiplicity $\{m_1,m_2,\ldots m_\ell\}$}
($m_i\in\bbz$) if the inclusion from $M_i$ to $W$ induces an
isomorphism on $H_1(M_i;\bbq)$ such that a generator 1 in
$H_1(M_i)/torsion$ is sent to $m_i$ in $H_1(W)/torsion$ and $W$
admits a rational $n$-Lagrangian with rational $n$-duals. In
addition to this, if the rational $n$-Lagrangian is the image of a
rational $(n+1)$-Lagrangian under the covering map, then $W$ is
called a \emph{rational $(n.5)$-cylinder of multiplicity
$\{m_1,m_2,\ldots,m_\ell\}$}.
\end{enumerate}
\end{defn}

\begin{rem}
\label{rem:n-cylinder}
\begin{enumerate}
\item By naturality of intersection forms and commutativity of the
diagram below, the image of an $n$-Lagrangian in $H_2(W)/I$
becomes a metabolizer for the $\ov{\la_0}$ and it
follows that the signature of $W$ is zero. The same is true for a
rational $n$-Lagrangian.
$$
\begin{diagram}\dgARROWLENGTH=1.0em
\node{H_2(W^{(n)})} \arrow{s} \node{\times} \node{H_2(W^{(n)})}
\arrow{s} \node{\lra} \node[2]{\bbz[\1(W)/\1(W)^{(n)}]} \arrow{s}\\
\node{H_2(W)/I} \node{\times} \node{H_2(W)/I} \node{\lra}
\node[2]{\bbz}
\end{diagram}
$$

\item Using twisted local coefficients, $H_2(W^{(n)})$ is
identified with $H_2(W;\bbz[\1(W)/\1(W)^{(n)}])$ as (right)
$\bbz[\1(W)/\1(W)^{(n)}]$-modules.

\item An $n$-cylinder is, in particular, a rational $n$-cylinder
of multiplicity $\{1,1,\ldots,1\}$.

\item A spin homology cobordism, say $W$, between two closed,
connected, oriented 3-manifolds is an $n$-cylinder since
$H_2(W)/I=0$.
\end{enumerate}
\end{rem}

The notion of $n$-cylinder extends the notion of $n$-solution in
\cite[Definitions 1.2,8.5,8.7]{COT1}. This is clear from the
following proposition.

\begin{prop}
\label{prop:n-cylider and n-solution} A 4-manifold $W$ with one
boundary component is an $n$-cylinder if and only if it is an
$n$-solution.
\end{prop}
\begin{proof}
Suppose $V$ is a compact, connected, spin 4-manifold with
$\partial V = N$ where $N$ is connected, $H_1(N) \cong \bbz$, and
the inclusion from $N$ to $V$ induces an isomorphism on $H_1(N)$.
To prove the proposition, we only need to show $J \equiv
\Image\{\inc_* : H_2(N) \lra H_2(V)\} = 0$ where $\inc_*$ is the
homomorphism induced by the inclusion. Since the inclusion induces
an isomorphism on the first homology, $\inc^* : H^1(V) \lra
H^1(N)$ is an isomorphism. By Poincar$\acute{\text{e}}$ duality,
the boundary map $H_3(V,N) \lra H_2(N)$ is an isomorphism. Hence
$\inc_* : H_2(N) \lra H_2(V)$ is a zero homomorphism in the
homology long exact sequence of the pair $(V,N)$.
\end{proof}

\begin{rem}
\label{rem:rational n-cylinder and rational n-solution} In a
similar fashion, one can easily see that a 4-manifold $W$ with one
boundary component is a rational $n$-cylinder if and only if it is
a rational $n$-solution (see \cite[Definition 4.1]{COT1}).
\end{rem}

\noindent Using $n$-cylinders, we define $n$-solvequivalence
between 3-manifolds and between knots.

\begin{defn}
\label{defn:n-solvequivalence}Let $n\in \frac12\bbn_0$. Let $M_1$
and $M_2$ be closed, connected, oriented 3-manifolds.
\begin{enumerate}
\item If there
exists an $n$-cylinder $W$ such that $\partial W = M_1\coprod
-M_2$, then $M_1$ is \emph{$n$-solvequivalent to $M_2$ via $W$}.

\item If there exists a rational $n$-cylinder (of multiplicity
$\{m_1,-m_2\}$) $W$ such that $\partial W = M_1\coprod -M_2$, then
$M_1$ is \emph{rationally $n$-solvequivalent to $M_2$ via $W$ (of
multiplicity $\{m_1,m_2\}$)}.

\item For given two knots $K_i$ $(i=1,2)$, if the zero framed
surgery on $S^3$ along $K_1$ is (rationally) $n$-solvequivalent to
the zero framed surgery on $S^3$ along $K_2$ via $W$, then we say
$K_1$ is \emph{(rationally) $n$-solvequivalent to $K_2$ via $W$}.
\end{enumerate}
\end{defn}

\noindent It is not hard to show that $n$-solvequivalence is an
equivalence relation. Here we give a proof only for transitivity.

\begin{prop}
\label{prop:transitivity} Let $L$, $M$, and $N$ be closed,
connected, oriented 3-manifolds with $H_1(L) \cong H_1(M) \cong
H_1(N) \cong \bbz$. Let $n\in \frac12 \bbn_0$. Suppose $L$ is
$n$-solvequivalent to $M$ via $V$ and $M$ is $n$-solvequivalent to
$N$ via $W$. Then $L$ is $n$-solvequivalent to $N$ via $V\cup_M
W$.
\end{prop}
\begin{proof}
Let $n\in \bbn$. Let $X \equiv V\cup_M W$. Then $\partial X = L
\coprod -N$. We need to show that $X$ is an $n$-cylinder. The
proof for the condition on the first homology groups is not hard,
hence left for the readers. Denote the inclusion map from $A$ to
$B$ by $i_{A,B}$ for two topological spaces $A \subset B$. Using
the long exact sequences of pairs, one can prove that
$H_2(V)/(i_{\partial V,V})_*(H_2(\partial V)) \cong
H_2(V)/(i_{L,V})_*(H_2(L)) \cong H_2(V)/(i_{M,V})_*(H_2(M))$. We
have similar isomorphisms for $(W,M\coprod -N)$. Now from the
Mayer-Vietoris sequence for $V$ and $W$ along $M$, one can see
that
$$
\(H_2(V)/(i_{\partial V, V})_*(H_2(\partial V))\) \op
\(H_2(W)/(i_{\partial W,W})_*(H_2(\partial W))\) \cong
H_2(X)/(i_{\partial X,X})_*(H_2(\partial X)).
$$
Since $(i_{V,X})_*$ and $(i_{W,X})_*$ map $\1(V)^{(n)}$ and
$\1(W)^{(n)}$ into $\1(X)^{(n)}$, respectively, and intersection
forms are natural, the ``union" of the $n$-Lagrangian with
$n$-duals for $V$ and the $n$-Lagrangian with $n$-duals for $W$
constitutes an $n$-Lagrangian with $n$-duals for $X$.

The proof for the case when $n$ is a half-integer is left to the
reader.
\end{proof}

\noindent The notion of $n$-solvability defined by COT suggests an
equivalence relation on $\SC$ wherein $K_0 \thicksim K_1$ if and
only if $K_0\#(-K_1)$ is $n$-solvable. This relation has not been
given a name. The following proposition reveals the close
connection, for knots, between $n$-solvequivalence and the
equivalence relation arising form COT's $n$-solvability. We
initially expected that these two equivalence relations were in
fact identical. However, we have not proved this and now believe
that they may be slightly different.

\begin{prop}
\label{prop:old and new} Let $n\in \frac12\bbn_0$. For knots
$K_0$ and $K_1$, if $K_0\#(-K_1)$ is $n$-solvable, then $K_0$ is
$n$-solvequivalent to $K_1$.
\end{prop}
\begin{proof}
Suppose $n \in \bbn_0$. Let $M_0$, $M_1$, and $M$ denote the zero
surgeries on $S^3$ along $K_0$, $K_1$, and $K_0\#(-K_1)$,
respectively. We construct a (standard) cobordism between
$M_0\coprod (-M_1)$ and $M$. Take a product $(M_0\coprod -M_1)
\times [0,1]$. By attaching a 1-handle along $(M_0\coprod
-M_1)\times \{1\}$ we get a 4-manifold whose upper boundary is
$M_0\# (-M_1)$. Note that $M_0\# (-M_1)$ is also obtained by
taking zero framed surgery on $S^3$ along a split link $K_0\coprod
-K_1$. Next add a zero framed 2-handle along the upper boundary of
this 4-manifold such that the attaching map is an unknotted circle
which links $K_0$ and $-K_1$ once. Then the resulting 4-manifold
has upper boundary $M$ as may be seen by sliding the 2-handle
represented by zero surgery on $K_0$ over that of $K_1$. This
4-manifold is the cobordism $C$ between $M_0\coprod -M_1$ and $M$.
Since $K_0\#(-K_1)$ is $n$-solvable, $M$ bounds an $n$-solution
$V$. Let $W \equiv C\cup_M V$. We claim that $W$ is an
$n$-cylinder with $\partial W = M_0 \coprod -M_1$. This will
complete the proof.

Note that $\partial W = M_0\coprod -M_1$. One sees that $H_2(M_0)$
and $H_2(M_1)$ are generated by capped-off Seifert surfaces for
$K_0$ and $K_1$, respectively. It is easy to show that $H_2(C) =
H_2(M_0) \oplus H_2(M_1)$ $(\cong \bbz \oplus \bbz)$ and $H_1(C)
\cong H_1(M_1) \cong H_1(M_0) \cong H_1(M)$ $(\cong \bbz)$ using
Mayer-Vietoris sequences. Also one sees that the inclusion from
$M$ to $C$ induces an injection $\inc_* : H_2(M) \lra H_2(C)$
where $H_2(M) \cong \bbz$ and $\inc_*$ sends 1 $(\in \bbz)$ to
(1,1) $(\in \bbz \oplus \bbz)$. Consider the following
Mayer-Vietoris sequence :
$$
\cdots \to H_2(M) \xrightarrow{f} H_2(C)\oplus H_2(V) \to H_2(W)
\xrightarrow{g} H_1(M) \xrightarrow{h} H_1(C)\oplus H_1(V) \to
\cdots .
$$
Since $h$ is an injection, $g$ is a zero map. Since $H_3(V,M) \to
H_2(M)$ is a dual map of $H^1(V) \to H^1(M)$ which is an
isomorphism, by the long exact sequence of homology groups of the
pair $(V,M)$, the inclusion induced homomorphism $H_2(M) \to
H_2(V)$ is a zero map. Thus the image of $f$ in $H_2(V)$ is zero
and by our previous observation the image of $f$ in $H_2(C)$ is
isomorphic to $\bbz$ generated by $(1,1)$ in $H_2(C)$. Therefore
$H_2(W) \cong H_2(V)\oplus \bbz$. Furthermore the last $\bbz$
summand on the right hand side is exactly $\Image \{\inc_* :
H_2(\partial W) \lra H_2(W)\}$, which is denoted by $I$. Hence
$H_2(W)/I \cong H_2(V)$. Since $V$ is a subspace of $W$,
$\1(V)^{(n)}$ is a subgroup of $\1(W)^{(n)}$. This implies that
$k$-surfaces in $V$ are also $k$-surfaces in $W$ for every integer
$k$. Now by naturality of (equivariant) intersections forms, one
can prove that an $n$-Lagrangian with $n$-duals for $V$ maps to an
$n$-Lagrangian with $n$-duals for $W$.

Clearly $W$ is a compact, connected 4-manifold, and it only
remains to show $W$ is spin. But this is obvious since $C$ and $V$
are spin and when we take the union of $C$ and $V$ we can adjust
spin structure of either of $C$ and $V$ to make $W$ spin.

Finally, in the case $K_0\# (-K_1)$ is $(n.5)$-solvable via $V$,
$K_0$ is $(n.5)$-solvequivalent to $K_1$ via $W$ where $W$ is
constructed as above. The argument for a proof for this goes the
same as above, and one just needs to notice that an
$(n+1)$-Lagrangian with $n$-duals for $V$ maps to an
$(n+1)$-Lagrangian with $n$-duals for $W$.
\end{proof}

\begin{rem}
\label{rem:old and new} There are a couple of natural arguments
for attempting to prove the converse of Proposition~\ref{prop:old
and new}, but they do not work completely since capped-off Seifert
surfaces for knots do not lift to higher covers of 4-manifolds in
the arguments.
\end{rem}

\section{Homology groups of an $n$-cylinder and its boundary}
\label{sec:homology} Let $\G$ be a poly-(torsion-free-abelian)
group (abbreviated PTFA). Then $\bbz\G$ is an Ore Domain and thus
embeds in its classical right ring of quotients $\SK_\G$, which is
a skew field \cite[Proposition 2.5]{COT1}. The skew field $\SK_\G$
is a $\bbz\G$-bimodule and has useful properties. In particular,
$\SK_\G$ is flat as a left $\bbz\G$-module (see \cite[Proposition
II.3.5]{Ste}), and every module over $\SK_\G$ is a free module
with a well defined rank over $\SK_\G$. The rank of any
$\bbz\G$-module $M$ is then defined to be the $\SK_\G$-rank of
$M\otimes_{\mathbb{Z}\G} \SK_\G$. Now we investigate $H_0$, $H_1$,
and $H_2$ of an $n$-cylinder $W$ and its boundary with
coefficients in $\bbz\G$ or $\SK_\G$. Throughout this section,
$\G$ denotes a PTFA group and $\SK_\G$ its (skew) quotient field
of fractions. The following two basic propositions are due to
\cite{COT1}.

\begin{prop}\cite[Proposition 2.9]{COT1}
\label{prop:H_0} Suppose $X$ is a $CW$-complex and there is a
homomorphism $\phi : \1(X) \lra \G$. Suppose $\psi : \bbz\G \to
\SR$ defines $\SR$ as a $\bbz\G$-bimodule and some element of the
augmentation ideal of $\bbz[\1(X)]$ is invertible in $\SR$. Then
$H_0(X;\SR)=0$. In particular, if $\phi : \1(X) \to \G$ is a
nontrivial coefficient system, then $H_0(X;\SK_\G) = 0$.
\end{prop}
\begin{prop}\cite[Proposition 2.11]{COT1}
\label{prop:H_1} Suppose $X$ is a $CW$-complex such that $\1(X)$
is finitely generated, and $\phi : \1(X) \lra \G$ is a nontrivial
coefficient system. Then
$$
\rk_{\SK_\G} H_1(X;\SK_\G)\le \b_1(X)-1.
$$
In particular, if $\b_1(X) = 1$, then $H_1(X;\SK_\G) = 0$. That
is, $H_1(X;\bbz\G)$ is a $\bbz\G$-torsion module.
\end{prop}
Since a rational $n$-cylinder $W$ has $\beta_1(W)=1$ we have:
\begin{cor}
\label{cor:homology} Suppose $W$ is a rational $n$-cylinder with
$\partial W = \coprod^\ell_{i=1}M_i$ where each $M_i$ is a
connected component of $\partial W$. If $\phi : \1(W) \lra \G$ is
a nontrivial coefficient system, then $H_0(W;\SK_\G) =
H_1(W;\SK_\G) = 0$ and $H_0(M_i;\SK_\G) = H_1(M_i;\SK_\G) = 0$ for
all $i$. Moreover, $H_2(M_i;\SK_\G)=0$ for all $i$.
\end{cor}
\begin{proof} A short technical argument shows that the restriction
of the coefficient system to each $\pi_1(M_i)$ is non-trivial (see
the proof of \cite[Proposition 2.11]{COT1}). The proof for $H_0$
and $H_1$ follows from the previous two propositions and the
definition of a rational $n$-cylinder. Notice that
$H_2(M_i;\SK_\G) \cong H^1(M_i;\SK_\G) \cong H_1(M_i;\SK_\G) = 0$
by Poincar$\acute{\text{e}}$ duality and the universal coefficient
theorem.
\end{proof}

\noindent The following proposition about $H_2$ plays an essential
role in showing that von Neumann $\rho$-invariants obstruct
$n$-solvequivalence. It extends Proposition 4.~3 of
\cite{COT1}. For the case that $W$ has one
boundary component, refer
to Proposition 4.~3 in \cite{COT1}. Let $I_\bbq \equiv \Image
\{\inc_* : H_2(\partial W;\bbq) \to H_2(W;\bbq)\}$.

\begin{prop}
\label{prop:H_2} Suppose $W$ is a compact, connected, oriented
$4$-manifold with $\partial W = M_1 \coprod M$ where $M =
\coprod^\ell_{i=2}M_i$ and $M_i$ are connected, $i=1,2,\ldots ,
\ell$ $(\ell\ge 2)$, and, for all $i$, $H_1(M_i;\mathbb{Q})\cong
H_1(W;\mathbb{Q})\cong\mathbb{Q}$, induced by inclusion. Suppose
$\phi : \1(W) \lra \G$ is a nontrivial (PTFA) coefficient system.
Then
$$
\rk_{\SK_\G}\!H_2(W;\SK_\G) = \rk_\bbq (H_2(W;\bbq)/I_\bbq) =
(\rk_\bbq H_2(W;\bbq)) - (\ell -1).
$$
Moreover, suppose there are 2-dimensional surfaces $S_j$ and
continuous maps $f_j : S_j \lra W$ $(j\in J)$ that are lifted to
$\tilde{f}_j \lra W_\G$ where $W_\G$ is a regular $\G$-cover
associated to $\phi$. If $\{[f_j]\,|\,j\in J\}$ is linearly
independent in $H_2(W,M;\bbq)$, then $\{[\tilde{f}_j]\,|\,j\in
J\}$ is $\bbq\G$-linearly independent in $H_2(W,M;\bbq\G)$.
\end{prop}
\begin{proof}
Since $W$ has nontrivial boundary, we can choose a 3-dimensional
$CW$-complex structure for $(W,\partial W)$. Let $C_*(W)$ be a
cellular chain complex associated to a chosen $CW$-complex
structure with coefficients in $\bbq$. Let $C_*(W_\G)$ be the
corresponding free $\kg$ chain complex of $W_\G$ that is freely
generated on cells of $W$. Let $b_i
\equiv \rk_\kg\! H_i(W;\kg)$ and $\b_i \equiv \rk_\bbq
H_i(W;\bbq)$ for $1\le i \le 4$. In this proof all chain complexes
and homology groups are with coefficients in $\bbq$ unless
specified otherwise.

By Corollary~\ref{cor:homology}, $b_0 = b_1 = 0$. By Poincar$\ee$
duality and the universal coefficient theorem, $H_3(W;\kg) \cong
H^1(W,\partial W;\kg) \cong H_1(W,\partial W;\kg)$. Using the long
exact sequence for the pair $(W,\partial W)$ with coefficients in
$\kg$, Proposition \ref{prop:H_0} and Proposition \ref{prop:H_1},
one sees that $H_1(W,\partial W;\kg)= 0$. Here we also need that
the coefficient system restricted to each $M_i$ is non-trivial as
mentioned previously (see proof of \cite[Proposition 2.11]{COT1}).
Hence $b_3 = 0$.

On the other hand, $H_3(W) \cong H^1(W,\partial W) \cong
H_1(W,\partial W)$. In the following long exact sequence
$$
\cdots \to H_1(\partial W) \xrightarrow{f} H_1(W) \to
H_1(W,\partial W) \to H_0(\partial W) \xrightarrow{g} H_0(W) \to
H_0(W,\partial W) \to 0,
$$
the maps $f$ and $g$ are surjective. Therefore $H_1(W,\partial W)
\cong \bbq^{\ell-1}$ and $\b_3 = \ell -1$. It is clear that $\b_0
= \b_1 = 1$. Since the Euler characteristics of $C_*(W)$ and
$C_*(W_\G)$ are equal, we deduce that $b_2 = \b_2 - (\ell -1)$.

For $\rk_\bbq (H_2(W)/I_\bbq)$, use the following long exact
sequence :
$$
H_2(\partial W) \to H_2(W) \to H_2(W,\partial W) \to H_1(\partial
W) \xrightarrow{f} H_1(W).
$$
Since $f$ is a surjection from $\bbq^\ell$ to $\bbq$, we have an
exact sequence as follows.
$$
0 \to H_2(W)/I_\bbq \to H_2(W,\partial W) \to \bbq^{\ell -1} \to
0.
$$
Thus $\rk_\bbq (H_2(W)/I_\bbq) = \b_2 - (\ell -1)$. This completes
the first part of the proof.

For the second part, let $X$ be the one point union of $S_j$ using
some base paths. Let $f : X \lra W$ and $\tilde{f} : X \lra W_\G$
be maps which restrict to $f_j$ and $\tilde{f}_j$, respectively.
By taking mapping cylinders, we may think of $X$ as a
(2-dimensional) subcomplex of $W$ and $C_*(X)$ as a subcomplex of
$C_*(W)$. If we denote by $X_\G$ the regular $\G$-cover associated
to $\phi\circ f_*$ (which can be thought of as a subcomplex of
$W_\G$), then $C_*(X_\G)$ is a subcomplex of $C_*(W_\G)$. In fact,
by our hypothesis $f_j$ lifts to $\tilde{f}_j$, and this implies
that $\phi\circ f_*$ is trivial on $\1(X)$. Hence $X_\G$ is a
trivial cover which consists of $\G$ copies of $X$. Consider the
following commutative diagram where each row is exact.
$$
\begin{diagram}\dgARROWLENGTH=1.5em
\node{H_3(W,M;\bbq\G)} \arrow{e}
   \arrow{s} \node{H_3(W,M\cup X;\bbq\G)}
   \arrow{e,t}{\tilde{\partial}} \arrow{s} \node{H_2(X;\bbq\G)}
   \arrow{e,t}{\tilde{f}_*} \arrow{s}
   \node{H_2(W,M;\bbq\G)}\arrow{s}\\
\node{H_3(W,M)} \arrow{e}
   \node{H_3(W,M\cup X)}
   \arrow{e,t}{\partial} \node{H_2(X)}
   \arrow{e,t}{f_*}
   \node{H_2(W,M)}
\end{diagram}
$$
Since $X_\G$ is a trivial cover, $H_2(X;\bbq\G)$ is a free
$\bbq\G$-module on $\{[S_j]\}$. Thus to complete the proof, we
need to show $\tilde{f}_*$ is injective. Since $H_3(W,M) \cong
H^1(W,M_1) \cong H_1(W,M_1) = 0$ and $C_*(W,M)$ is a 3-dimensional
chain complex, $\partial_\# : C_3(W,M) \lra C_2(W,M)$ is
injective. Notice that $C_*(W,M)$ can be identified with
$C_*(W_\G,M_\G)\otimes_{\bbq\G}\bbq$ (similarly for $C_*(X)$ and
$C_*(W,M\cup X)$). Then by \cite[pg.305]{Str} (or see
\cite[Proposition 2.4]{COT1}), $\partial_\# : C_3(W_\G,M_\G) \lra
C_2(W_\G,M_\G)$ is also injective. Hence we obtain
$H_3(W,M;\bbq\G)=0$. Now it suffices to show $H_3(W,M\cup X
;\bbq\G)=0$. Since $f_*$ is injective by our hypothesis,
$H_3(W,M\cup X) = 0$. Since $C_*(W,M\cup X)$ is a 3-dimensional
chain complex, this implies that $\partial_\# : C_3(W,M\cup X)
\lra C_2(W,M\cup X)$ is injective. Once again by
\cite[pg.305]{Str}, $\partial_\# : C_3(W_\G,M_\G\cup X_\G) \lra
C_2(W_\G,M_\G\cup X_\G)$ is injective. Therefore $H_3(W,M\cup
X;\bbq\G) = 0$.
\end{proof}

We now investigate the relationship between the first homology
groups of an $n$-cylinder and its boundary components. The
following lemma generalizes \cite[Lemma 4.5]{COT1}. It is the
linchpin in proving Theorem~\ref{thm:rank}, which establishes the
crucial connection between $n$-solvequivalence and homology.

\begin{lem}
\label{lem:exact} Let $W$ be a rational $n$-cylinder with $M$ as
one of it boundary components. Let $\G$ be an $(n-1)$-solvable
(PTFA) group and $\SR$ be a ring such that $\bbq\G \subset \SR
\subset \kg$. Suppose $\phi : \1(M) \lra \G$ is a nontrivial
coefficient system that extends to $\psi : \1(W) \lra \G$. Then
$$
TH_2(W,M;\SR) \xrightarrow{\partial} H_1(M;\SR) \xrightarrow{i_*}
H_1(W;\SR)
$$
is exact. (Here for an $\SR$-module $\SM$, $T\SM$ denotes the
$\SR$-torsion submodule of $\SM$.)
\end{lem}
\begin{proof}
We need to show that every element of $\Ker (i_*)$ is in the image of an
element of $TH_2(W,M;\SR)$. Let $m =
\frac12 \rk_\bbq(H_2(W;\bbq)/I_\bbq)$ where $I_\bbq \equiv \Image
\{\inc_* : H_2(\partial W;\bbq) \to H_2(W;\bbq)\}$. By
Proposition~\ref{prop:H_2}, $\rk_\kg H_2(W;\SR) = 2m$. Let
$\{\ell_1,\ell_2,\ldots, \ell_m\}$ generate a rational
$n$-Lagrangian for $W$ and $\{d_1,d_2,\ldots, d_m\}$ be its
$n$-duals. Since $\G$ is $(n-1)$-solvable, $\psi$ descends to
$\psi' : \1(W)/\1(W)^{(n)} \lra \G$. We denote by $\ell_i'$ and
$d_i'$ the images of $\ell_i$ and $d_i$ in $H_2(W;\SR)$. By
naturality of intersection forms, the intersection form $\la$
defined on $H_2(W;\SR)$ vanishes on the module generated by
$\{\ell_1',\ell_2',\ldots, \ell_m'\}$. Let $\SR^m\op \SR^m$ be the
free module on $\{\ell_i',d_i'\}$.
The following composition
$$
\SR^m\op \SR^m \xrightarrow{j_*} H_2(W;\SR) \xrightarrow{\la}
H_2(W;\SR)^* \xrightarrow{j^*} (\SR^m\op \SR^m)^*
$$
is represented by a block matrix
$$
\left(\begin{matrix} 0 & I\cr I & X
\end{matrix}\right).
$$
This matrix has an inverse which is
$$
\left(\begin{matrix} -X & I\cr I & 0
\end{matrix}\right).
$$
Thus the composition is an isomorphism. This implies that $j_*$ is
a monomorphism and $j^*$ is a (split) epimorphism. Since $j^*$ is
a split epimorphism between the free $\SR$-modules of the same
$\kg$-rank, and $\SR$ is an integral domain, $j^*$ is an
isomorphism. Hence $\la$ is a surjection. Now consider the
following diagram where the row is the exact sequence of the pair
$(W,M)$.
$$
\begin{diagram}\dgARROWLENGTH=1.2em
\node[4]{H_2(\partial W,M;\SR)} \arrow{s,r}{f_*}\\
\node{H_2(W;\SR)} \arrow[3]{e,t}{\pi_*} \arrow[3]{se,r}{\la}
\node[3]{H_2(W,M;\SR)}
\arrow[3]{e,t}{\partial} \arrow{s,r}{g_*} \node[3]{H_1(M;\SR)}
\arrow[3]{e,t}{i_*} \node[3]{H_1(W;\SR)}\\
\node[4]{H_2(W,\partial W;\SR)} \arrow{s,r}{\PD}\\
\node[4]{H^2(W;\SR)} \arrow{s,r}{\k}\\
\node[4]{H_2(W;\SR)^*}.
\end{diagram}
$$
 We claim that $\Ker (\k\circ \PD \circ g_*)$ is $\SR$-torsion.
$\Ker(\k)$ is $\SR$-torsion since it is a split surjection between
$\SR$-modules of the same rank over $\kg$. $\PD$ is an isomorphism
by Poincar$\ee$ duality. $\Ker(g_*) = \Image (f_*)$ and
$H_2(\partial W, M;\SR) = H_2(M';\SR)$ where $M'$ is the disjoint
union of the boundary components of $W$ except for $M$. Since
$H_2(M';\SR)$ is $\SR$-torsion by the flatness of $\kg$ over $\SR$
and Corollary~\ref{cor:homology} (once again we need to know that
the restricted coefficient system is non-trivial), $\Ker (g_*)$ is
$\SR$-torsion. Combining these, one can deduce that $\Ker (\k\circ
\PD \circ g_*)$ is $\SR$-torsion.

Suppose $p\in \Ker(i_*) \subset H_1(M;\SR)$. Then there exists $x
\in H_2(W,M;\SR)$ such that $\partial(x) = p$. Let $y\in
\la^{-1}((\k\circ\PD\circ g_*)(x))$. Then $x-\pi_*(y) \in \Ker(\k
\circ \PD \circ g_*)$. Hence $x-\pi_*(y)$ is $\SR$-torsion and
$\partial(x-\pi_*(y)) = \partial(x) = p$.
\end{proof}

Under the same hypotheses as in Lemma~\ref{lem:exact}, there
exists a non-singular linking form $\SB\ell : H_1(M;\SR) \lra
H_1(M;\SR)^\# \equiv \overline{\Hom_\SR(H_1(M;\SR), \kg/\SR)}$ by
\cite[Theorem 2.13]{COT1}. This definition, and a proof of
nonsingularity will be included in our proof of
Proposition\ref{prop:P and P perb}. For an $\SR$-submodule $P$ of
$H_1(M;\SR)$, we define $P^\bot \equiv \{x\in H_1(M;\SR)\,\,|\,\,
\SB\ell(x)(y) = 0, \forall y \in P\}$, clearly an $\SR$-submodule.

\begin{prop}
\label{prop:P and P perb} Suppose the same hypotheses as in
Lemma~\ref{lem:exact}. Suppose $\SR$ is a PID. Then for $P\equiv
\Ker \{i_* : H_1(M;\SR) \lra H_1(W;\SR)\}$, $P\subset P^\bot$.
Moreover, if $M = \partial W$, then $ P = P^\bot$.
\end{prop}
\begin{proof}
Consider the following commutative diagram.
$$
\begin{diagram}\dgARROWLENGTH=1.2em
\node{TH_2(W,M;\SR)} \arrow{e,t}{\partial} \arrow{s,r}{g_*}
\node{H_1(M;\SR)} \arrow{e,t}{i_*} \arrow[2]{s,r}{\PD}
\node{H_1(W;\SR)}\\
\node{TH_2(W,\partial W;\SR)} \arrow{s,r}{\PD}\\
\node{TH^2(W;\SR)} \arrow{e,t}{i^*} \arrow{s,r}{B^{-1}}
\node{H^2(M;\SR)} \arrow{s,r}{B^{-1}}\\
\node{H^1(W;\kg/\SR)} \arrow{e,t}{j^*} \arrow{s,r}{\k}
\node{H^1(M;\kg/\SR)} \arrow{s,r}{\k}\\
\node{H_1(W;\SR)^\#} \arrow{e,t}{i^\#} \node{H_1(M;\SR)^\#}
\end{diagram}
$$
Let $\br : TH_2(W,M;\SR) \lra H_1(W;\SR)^\#$ be the composition of
the maps on the left column and $B\ell : H_1(M;\SR) \lra
H_1(M;\SR)^\#$ the composition of the maps on the right column.
$\PD$ are induced by Poincar$\ee$ duality (hence isomorphisms),
and $\k$ are the Kronecker evaluation maps. $B^{-1}$ are the
inverses of the Bockstein homomorphisms. For the existence of
$B^{-1}$, we need to show that the Bockstein homomorphisms are
isomorphisms. First, we consider the following exact sequence.
$$
H^1(W;\kg) \to H^1(W;\kg/\SR) \xrightarrow{B} H^2(W;\SR) \to
H^2(W;\kg).
$$
Notice that $H^1(W;\kg) \cong \Hom_\kg(H_1(W;\kg),\kg) = 0$ by
Corollary~\ref{cor:homology}. Since $H^2(W;\kg)$ is $\SR$-torsion
free and $H^1(W;\kg/\SR)$ is $\SR$-torsion, one sees that $B$ is
an isomorphism onto $TH^2(W;\SR)$. Secondly, one can prove that
$H^1(M;\kg) = H^2(M;\kg) = 0$ using Corollary~\ref{cor:homology}.
Hence $B : H^1(M;\kg/\SR) \lra H^2(M;\SR)$ is an isomorphism in
the following long exact sequence.
$$
H^1(M;\kg) \to H^1(M;\kg/\SR) \xrightarrow{B} H^2(M;\SR) \to
H^2(M;\kg).
$$

If $x\in P$, then $x = \partial(y)$ for some $y \in TH_2(W,M;\SR)$
by Lemma~\ref{lem:exact}. Thus $B\ell(x) = (i^\#\circ\br)(y)$. For
every $p\in P$, $B\ell(x)(p) = (i^\#\br)(y)(p) = \br(y)(i_*(p)) =
\br(y)(0) = 0$. Therefore $x\in P^\bot$, and $P\subset P^\bot$.
The proof for the case when $M = \partial W$ follows from
\cite[Theorem 4.4]{COT1}.
\end{proof}

\begin{rem}
\label{rem:rank} If $W$ has more than one boundary component, then
in general $P \ne P^\bot$. For example, if $W= M\times [0,1]$, $P
= 0$ and $P^\bot  = H_1(M;\SR)$.
\end{rem}

Suppose $\G$ is a PTFA group with $H_1(\G) \cong \G/[\G,\G] \cong
\bbz$. Then its commutator subgroup, $[\G,\G]$, is also PTFA and
$\bbz[\G,\G]$ embeds into its (skew) quotient field of fractions
denoted by $\bbk$. Then we have a PID $\bbk[t^{\pm 1}]$ such that
$\bbz\G \subset \bbk[t^{\pm 1}] \subset \kg$ where $t$ is
identified with the generator of $\G/[\G,\G]$.

The following generalizes \cite[Theorem 6.4]{CT}. The proof is the
same once equipped with Proposition~\ref{prop:P and P perb}.
\begin{thm}
\label{thm:rank} Let $M$ be zero surgery on a knot $K$ in $S^3$.
Let $W$ be an $n$-cylinder with $M$ as one of its boundary
components. Suppose $\G$ is an $(n-1)$-solvable PTFA group with
$H_1(\G) \cong \bbz$. Suppose $\phi : \1(W) \lra \G$ induces an
isomorphism upon abelianization. Let $d \equiv \rk_\bbq
H_1(M_\infty;\bbq)$ where $M_\infty$ is the (universal) infinite
cyclic cover of $M$. Then
$$
\rk_\bbk \Image\{i_* : H_1(M;\bbk[t^{\pm 1}]) \lra
H_1(W;\bbk[t^{\pm 1}])\} \ge (d-2)/2 \phantom{a} \text{if}
\phantom{a} n>1
$$
and this rank is at least $d/2$ if $n=1$.
\end{thm}
\begin{proof}
Let $\SR \equiv \bbk[t^{\pm 1}]$. By \cite[Theorem 2.13]{COT1},
there exists a non-singular Blanchfield linking form $B\ell :
H_1(M;\SR) \lra H_1(M;\SR)^\# \equiv
\overline{\Hom_\SR(H_1(M;\SR),\kg/\SR)}$. Let $P\equiv \Ker \{i_*
: H_1(M;\SR) \lra H_1(W;\SR)\}$, $Q \equiv \Image\{i_* :
H_1(M;\SR) \lra H_1(W;\SR)\}$, and $\SA \equiv H_1(M;\SR)$. By
Proposition~\ref{prop:P and P perb}, $P\subset P^\bot$ with
respect to the Blanchfield linking form. This gives us a
well-defined map $f : P \lra (\SA/P)^\#$ induced from the
Blanchfield linking form. Since the Blanchfield linking form is
non-singular, $f$ is a monomorphism. Hence $\rk_\bbk P \le
\rk_\bbk (\SA/P)^\#$.

We claim that $\rk_\bbk \SM = \rk_\bbk (\SM)^\#$ for every
finitely generated $\SR$-module $\SM$. (Here $(\SM)^\# \equiv
\overline{\Hom_\SR(\SM,\kg/\SR)}$.) Since $\SR$ is a PID, it is
enough to show this for the case when $\SM$ is a cyclic
$\SR$-module, i.~e.~, $\SM = \SR/p(t)$ where $p(t) \in \SR$. For
this cyclic case, $(\SR/p(t))^\# \cong \SR/\overline{p}(t)$ where
$\overline{p}(t)$ is obtained by taking involution of $p(t)$.
Since $\rk_\bbk \SR/p(t)$ is the degree of $p(t)$ for every $p(t)
\in \SR$ and $p(t)$ and $\overline{p}(t)$ have the same degree,
$\rk_\bbk \SR/p(t) = \rk_\bbk \SR/\overline{p}(t)$. Using this we
can deduce that
$$
\rk_\bbk P \le \rk_\bbk (\SA/P)^\# = \rk_\bbk \SA/P = \rk_\bbk \SA
- \rk_\bbk P
$$
Hence $\rk_\bbk P \le \frac12 \rk_\bbk\SA$. Since $Q \cong \SA/P$
as $\SR$-modules, we have
$$
\rk_\bbk Q = \rk_\bbk \SA - \rk_\bbk P \ge
\rk_\bbk \SA - \frac12 \rk_\bbk \SA = \frac12 \rk_\bbk \SA.
$$
Thus we only need to show $\rk_\bbk \SA \ge d-2$ if $n>1$ and
$\rk_\bbk \SA = d$ if $n=1$.

If $n=1$, $\G \cong \bbz$, $\bbk = \bbq$, and $\SR = \bbq[t^{\pm
1}]$. In this case $\SA$ is the rational Alexander module, which
is $H_1(M_\infty;\bbq)$. Hence $\rk_\bbk \SA = d$. Suppose $n>1$.
Observe that $\SA$ is obtained from $H_1(\sk;\SR)$ by killing the
$\SR$-submodule generated by the longitude, say $\ell$. By
\cite[Corollary 4.8]{C}, $\rk_\bbk H_1(\sk;\SR) \ge d-1$. Since
$(t-1)_*\ell = 0$ in $H_1(\sk;\SR)$, the submodule generated by
$\ell$ is isomorphic with $\bbk[t^{\pm 1}]/(t-1) \cong \bbk$.
Hence $\rk_\bbk \SA \ge d-2$.
\end{proof}

\section{Obstructions to $n$-solvequivalence}
\label{sec:obstruction} We use the information about twisted
homology groups obtained in Section~\ref{sec:homology} to get
obstructions to $n$-solvequivalence. Our obstructions will be
vanishing of von Neumann $\rho$-invariants. The von Neumann
$\rho$-invariants were firstly used by Cochran-Orr-Teichner to
give obstructions to $n$-solvability (see Theorem $4.2$ and $4.6$
in \cite{COT1}). We begin by giving a very brief explanation about
von Neumann $\rho$-invariants. For more details on von Neumann
$\rho$-invariants, the readers are referred to Section 5 in
\cite{COT1} and Section 2 in \cite{CT}. Let $M$ be a compact,
oriented 3-manifold. If we have a representation $\phi : \1(M)
\lra \G$ for a group $\G$, then the \emph{von Neumann
$\rho$-invariant} or \emph{reduced $L^{(2)}$-signature}
$\rho(M,\phi)$ ($\in \bbr$) is defined. If $(M,\phi) = \partial
(W,\psi)$ for some compact, oriented 4-manifold $W$ and a
homomorphism $\psi : \1(W) \lra \G$, then we have $\rho(M,\phi) =
\s^{(2)}_\G(W,\psi) - \s_0(W)$ where $\s^{(2)}_\G(W,\psi)$ is the
$L^{(2)}$-signature of the intersection form defined on
$H_2(W;\bbz\G)$ twisted by $\psi$ and $\s_0(W)$ is the ordinary
signature of $W$. If $\G$ is a PTFA group, then as we have seen
before $\G$ embeds into the (skew) quotient field of fractions
$\kg$ and $\s^{(2)}_\G$ can be thought of as a homomorphism from
$L^0(\kg)$ to $\bbr$. Some useful properties of von Neumann
$\rho$-invariants due to COT are given below. One can find
detailed proofs or explanations in \cite[Section 5]{COT1}.

\begin{prop}
\label{prop:rho invariants}Let $M$ be as above and $\G$ a PTFA
group. Suppose we have a homomorphism $\phi : \1(M) \lra \G$.
\begin{itemize}
\item [(1)] If $(M,\phi) = \partial (W,\psi)$ for some compact,
spin 4-manifold $W$ and $H_2(W;\kg)$ has a half-rank summand on
which the (equivariant) intersection form vanishes, then
$\rho(M,\phi) = 0$. In fact, $\s^{(2)}_\G(W,\psi) = \s_0(W) = 0$.

\item [(2)] If $\phi$ factors through $\phi' : \1(M) \lra \G'$ where
$\G'$ is a subgroup of $\G$, then $\rho(M,\phi') = \rho(M,\phi)$.

\item [(3)] If $\phi$ is trivial, then $\rho(M,\phi) = 0$.

\item [(4)] If $\G=\bbz$ and $\phi$ extends over $W$ nontrivially
for some compact, spin 4-manifold $W$, then $\rho(M,\phi) =
\int_{\om \in S^1} \s(h(\om)) \mathrm{d}\om - \s_0(W)$ where $h$
is the matrix representing the intersection form on
$H_2(W;\bbc[t,t^{-1}])/\text{(torsion)}$.
\end{itemize}
\end{prop}

\noindent Now we are ready to give obstructions to
$n$-solvequivalence. In fact, we give a more general theorem.
\begin{thm}
\label{thm:rho=0} Let $M_i$ ($i=1,2,\ldots ,\ell)$ be closed,
connected, oriented 3-manifolds with $H_1(M_i;\bbq) \cong \bbq$
for all $i$. Let $\G$ be an $n$-solvable group. Suppose there
exists a rational $(n.5)$-cylinder $W$ with $\partial W =
\coprod^\ell_{i=1} M_i$ (of arbitrary multiplicity) and we have
representations $\phi_i : \1(M_i) \lra \G$. If $\phi_i$ extends to
the same homomorphism $\phi : \1(W) \lra \G$ for all $i$, then
$$
\sum^\ell_{i=1}\rho(M_i,\phi_i) = 0.
$$
\end{thm}

\noindent The following is an easy corollary of
Theorem~\ref{thm:rho=0}, hence the proof is omitted.
\begin{cor}
\label{cor:rho=0} Suppose $K_0$ is $(n.5)$-solvequivalent to $K_1$
via $W$ and $\G$ is an $n$-solvable group. Suppose we have
representations $\phi_i : \1(M_i) \lra \G$ for $i=0,1$ where $M_i$
is zero surgery on $S^3$ along $K_i$, $i=0,1$. If $\phi_i$ extends
to the same homomorphism $\phi : \1(W) \lra \G$ for $i=1,2$, then
$$
\rho(M_0,\phi_0) = \rho(M_1,\phi_1).
$$
\end{cor}

\begin{proof}[Proof of Theorem~\ref{thm:rho=0}]
Since $W$ is spin, $\s_0(W) = 0$. Since
$\sum^\ell_{i=1}\rho(M_i,\phi_i) = \s^{(2)}_\G(W,\phi) -\s_0(W)$,
we need to show $\s^{(2)}_\G(W,\phi)=0$. Let $\rk_\bbq
(H_2(W;\bbq)/I_\bbq) = 2m$ where $I_\bbq \equiv \Image\{\inc_* :
H_2(\partial W;\bbq) \lra H_2(W;\bbq)\}$. Let $L$ be a rational
$(n+1)$-Lagrangian for $W$ generated by $(n+1)$-surfaces
$\{\ell_1,\ell_2,\ldots ,\ell_m\}$. Since $\G$ is $n$-solvable,
$\phi$ factors through $\1(W)/\1(W)^{(n+1)}$. Thus if we let $\kg$
denote the (skew) quotient field of fractions of $\bbq\G$, then we
can take the image of $L$ in $H_2(W;\kg)$, which is denoted by
$L'$ (generated by $\{\ell'_1,\ell'_2,\ldots ,\ell'_m\}$). By
naturality of intersection form, the intersection form $\la' :
H_2(W;\kg)\times H_2(W;\kg) \to \kg$ vanishes on $L'$. By
Proposition~\ref{prop:H_2}, $\rk_\kg\! H_2(W;\kg) = 2m$. Therefore
if we show that $\{\ell'_1,\ell'_2,\ldots ,\ell'_m\}$ is linearly
independent in $H_2(W;\kg)$, then $\rk_\kg L' = m = \frac12\rk_\kg
H_2(W;\kg)$, which implies that $\s^{(2)}_\G(W,\phi) = 0$ by
Proposition~\ref{prop:rho invariants}(1).

For convenience we let $M\equiv \coprod^\ell_{i=2} M_i$ and $S
\equiv \{\ell_1,\ell_2,\ldots ,\ell_m\}$. Since $S$ generates a
rational $(n+1)$-Lagrangian and $\#(S) = m = \frac12\rk_\bbq
(H_2(W;\bbq)/I_\bbq)$, the image of $S$ in $H_2(W;\bbq)/I_\bbq$ is
linearly independent. By investigating the long exact sequences of
homology groups of the pairs $(W,\partial W)$ and $(W,M)$, one
easily sees that $I_\bbq = \Image\{\inc_* : H_2(M;\bbq) \lra
H_2(W;\bbq)\}$. It follows that $H_2(W;\bbq)/I_\bbq$ can be
identified with a $\bbq$-subspace of $H_2(W,M;\bbq)$. Hence the
image of $S$ in $H_2(W,M;\bbq)$ is linearly independent. By the
second part of Proposition~\ref{prop:H_2}, the image of $S$ in
$H_2(W,M;\bbq\G)$ is linearly independent. Hence the image of $S$
in $H_2(W,M;\kg)$ is also linearly independent. By
Corollary~\ref{cor:homology}, $H_1(M;\kg) = H_2(M;\kg) = 0$. This
implies that $H_2(W;\kg) \cong H_2(W,M;\kg)$. Therefore
$\{\ell'_1,\ell'_2,\ldots ,\ell'_m\}$ is linearly independent in
$H_2(W;\kg)$.
\end{proof}

\section{Main Theorem}
\label{sec:main} In this section we prove the main theorem, that
simultaneously generalizes the main theorems of \cite{CT} and
\cite{L}\cite{K}. The strength of our theorem lies in being able
to work within the class of knots that have the same fixed
classical Seifert matrix (and Alexander module) and the same
\emph{higher-order} analogues of these. Therefore, after stating
the theorem, we review the higher-order Alexander modules. Then we
state and prove two theorems that are used at the end of the
section to prove the main theorem. The first of these,
Theorem~\ref{thm:injectivity}, is an important technical result
that significantly generalizes  \cite[Theorem 4.3]{CT}. Due to its
difficult technical nature, the proof of our
Theorem~\ref{thm:injectivity} is not completed until the next
section.

\begin{thm}[{\bf Main Theorem}]
\label{thm:main} Let $n\in\bbn$. Let $K$ be a knot whose Alexander
polynomial has degree greater than 2 (if $n=1$ then degree equal to
$2$ is allowed). Then there is an infinite family of knots $\{K_i
| i\in \bbn_0\}$ with $K_0 = K$ such that:

\begin{itemize}
\item [(1)] For each $i$, $K_i - K$ is $n$-solvable. In
particular, $K_i$ is $n$-solvequivalent to $K$. Moreover, $K_i$
and $K$ cobound, in $S^3\times\[0,1\]$, a smoothly embedded symmetric
Grope of height $n+2$.

\item [(2)] If $i\ne j$, $K_i$ is not $(n.5)$-solvequivalent to
$K_j$. In particular, $K_i-K_j$ is not $(n.5)$-solvable, and $K_i$
and $K_j$ do not cobound, in $S^3\times\[0,1\]$, any embedded
symmetric Grope of height $(n+2.5)$.

\item [(3)] Each $K_i$ has the same $m$-th integral higher-order
Alexander module as $K$, for $m=0,1,\ldots, n-1$. Indeed, if $G_i$
and $G$ denote the knot groups of $K_i$ and $K$ respectively, then
there is an isomorphism $G_i/(G_i)^{(n+1)}\to G/G^{(n+1)}$ that
preserves the peripheral structures.

\item [(4)] Each $K_i$ has the same \emph{$m^{th}$-order Seifert
presentation} as $K$,
 for $m=0,1,\ldots, n-1$.

\item [(5)] For each $i>j$, $K_i - K_j$ is of infinite order in
$\SF_{(n)}/\SF_{(n.5)}$.

\item [(6)] If $i,j>0$, $s(K_i)=s(K_j)$ and $\tau(K_i)=\tau(K_j)$,
where $s$ is Rasmussen's smooth concordance invariant and $\tau$
is the smooth concordance invariant of P.~Ozsv{\'a}th and
Z.~Szab{\'o} \cite{Ra}\cite{OS}.

\end{itemize}
\end{thm}

\noindent Since there certainly exist knots whose Alexander
polynomial has degree $4$, from ($5$) above we recover the result
of Cochran and Teichner.
\begin{cor}
\label{cor:cot filtration} \cite{CT} For every positive integer $n$,
$\SF_{(n)}/\SF_{(n.5)}$ has positive rank.
\end{cor}

\begin{cor} Let $n$ be a positive integer.
\label{cor:gropefiltrations}

\begin{itemize}
\item [1.] $\SG_{n+2}^{smooth}/\SG_{n+2.5}^{smooth}$ has
elements of infinite order represented by knots $K$ with $s(K)=0$
and $\tau(K)=0$.
\item [2.] The kernel of the homomorphism $ST:\SG_{n+2}^{smooth}\to
\mathbb{Z}\times \mathbb{Z}$ given by $ST(K)=(s(K),\tau(K))$ is
infinite.
\end{itemize}
\end{cor}

\begin{proof}[Proof of Corollary~\ref{cor:gropefiltrations}]
For fixed $n$, consider the knot $K_2-K_1$, that is $K_2\#(-K_1)$.
Note that, by part $(6)$ above, $\tau(K_2-K_1)=s(K_2-K_1)=0$. By
part $(1)$ above (applied first to $K_1$ and $K$ and then to $K_2$
and $K$), $K_1$ and $K_2$ cobound, in $S^3\times\[0,1\]$, a
smoothly embedded symmetric Grope of height $n+2$. It follows that
$K_1\#(-K_2)$ bounds a smoothly embedded symmetric Grope of height
$n+2$ in $B^4$. This is seen by choosing an embedded arc from
$K_1$ to $K_2$ in $S^3\times\[0,1\]$, contained in the first-stage
surface of the Grope and avoiding the boundaries of the
higher-stage surfaces. Deleting a tubular neighborhood of this
arc, one gets the desired Grope in $B^4$. But if some multiple of
$K_2\#(-K_1)$ were the boundary of a smoothly (or even
topologically) embedded symmetric Grope of height $n+2.5$ in
$B^4$, then by \cite[Theorem 8.11]{COT1}, it would be of finite
order in $\SF_{(n)}/\SF_{(n.5)}$ contradicting part $(5)$ above.
This completes the proof of the first part of the corollary.

For the second part of the corollary, just consider multiples of
$K_2-K_1$. Alternatively, note that the analysis above applies
equally well to any of the infinite set of knots $\{K_i-K_1$,
$i\geq 1\}$.
\end{proof}

\noindent The higher-order Alexander modules of a knot were
defined by the first author, and we give brief explanations here.
For more details refer to \cite{C}.

\begin{defn}\cite{C}
\label{defn:alexander module} Let $n$ be a nonnegative integer.
Let $K$ be a knot and $G \equiv \1(S^3\setminus K)$. The
\emph{$n$-th (integral) higher-order Alexander module}, $\Anz(K)$,
of a knot $K$ is the first (integral) homology group of the
covering space of $S^3\setminus K$ corresponding to $G^{(n+1)}$,
considered as a right $\bbz[G/G^{(n+1)}]$-module, i.~e.~,
$G^{(n+1)}/G^{(n+2)}$ as a right module over $\bbz[G/G^{(n+1)}]$.
\end{defn}

\noindent Using the local coefficient system $\pi:G\to
G/G^{(n+1)}$, one can see that $\Anz(K) \cong H_1(S^3\setminus
K;\bbz[G/G^{(n+1)}])$ as right $\bbz[G/G^{(n+1)}]$-modules. Notice
that $\SA_0^\bbz(K)$ is the classical Alexander module of $K$. Let
$\G_n \equiv G/G^{(n+1)}$. Let $\G'_n$ denote the commutator
subgroup of $\G_n$ (that is $G^{(1)}/G^{(n+1)})$. Note that when
$n=0$, $\G'_n=\{e\}$. Let $\Sigma$ be a Seifert surface for $K$, let
$E_K \equiv \sk$ and let $Y \equiv E_K - (\Sigma \times (-1,1))$. We
denote the two inclusions $\Sigma \to \Sigma\times\{\pm 1\} \to
\partial Y \subset Y$ by $i_+$ and $i_-$. Notice that the coefficient
system $\pi$, when restricted to $\pi_1(\Sigma)$ or $\pi_1(Y)$, has
image in $\G_n'$. Thus $\Sigma$ and $Y$ have naturally induced
$\mathbb{Z}\G_n'$-local coefficient systems and so we have the
inclusion-induced maps $(i_{\pm})_*:H_1(\Sigma;\mathbb{Z}\G'_n)\to
H_1(Y;\mathbb{Z}\G'_n)$. In the classical case $n=0$ these are
just the usual maps $(i_{\pm})_*:H_1(\Sigma;\mathbb{Z})\to
H_1(Y;\mathbb{Z})$ that determine the Seifert matrix and Seifert
form. For $n>0$, no analogue of the Seifert form has been found,
so these maps, though quite complicated, play a crucial role. In
the absence of a true analogue of a higher-order Seifert
\emph{form} they contain the relevant (integral) information
necessary to reconstruct the Alexander modules (but, as in the
classical case, contain \emph{more} information that the modules
alone). Therefore we make the following definition.

\begin{defn}
\label{defn:Seifert presentation} An \emph{$n^{th}$-order Seifert
presentation} of a knot $K$ is the ordered pair
$((i_+)_*,(i_-)_*)$ of maps
$(i_{\pm})_*:H_1(\Sigma;\mathbb{Z}\G'_n)\to
H_1(Y;\mathbb{Z}\G'_n)$ induced by the inclusion maps
$(i_{\pm}):\Sigma\to Y$ for some choice of Seifert surface
$\Sigma$ for $K$.
\end{defn}

Therefore, the case $n=1$ of our main theorem precisely recovers
that of \cite{Li} and the second author \cite{K}. In general,
since $H_1(\Sigma;\mathbb{Z}\G'_n)$ and $H_1(Y;\mathbb{Z}\G'_n)$
may not be free modules (or even finitely-generated), we cannot
speak of a Seifert \emph{matrix}. However after some localization
this situation can be remedied.

Recall from \cite{COT1}\cite{C} that there is a localized version
of $\Anz(K)$. We set $S_n \equiv \bbz\G_n' - \{0\}$, $n\ge 0$.
Then $S_n$ is a right divisor set of $\bbz\G_n$ and we let $\SR_n
\equiv (\bbz\G_n)(S_n)^{-1}$. It also follows that $\bbz\G_n'$ is
itself an Ore domain and embeds into its classical (skew) quotient
field of fractions, which we denote by $\bbk_n$. By
\cite[Proposition 3.2]{COT1}, $\SR_n$ is canonically identified
with the (skew) Laurent polynomial ring $\bbk_n[t^\pm]$, which is
a PID. In the classical case, $S_0=\mathbb{Z}-\{0\}$,
$\bbk_0=\mathbb{Q}$ and $\SR_0$ is $\mathbb{Q}[t^{\pm 1}]$. In
general, Cochran defines the localized higher-order Alexander
module \cite{C}.

\begin{defn}
\label{defn:localized Alexander module} The \emph{$n$-th localized
Alexander module} of a knot $K$ is $\SA_n(K) \equiv
H_1(\sk;\SR_n)$.
\end{defn}

Therefore the case $n=0$ gives the classical (rational) Alexander
module. Moreover the $n^{th}$-order Seifert presentation maps
$(i_{\pm})_*:H_1(\Sigma;\mathbb{Z}\G'_n)\to H_1(Y;\mathbb{Z}\G'_n)$
determine their localized counterparts
$(i_{\pm})_*:H_1(\Sigma;\bbk_n)\to H_1(Y;\bbk_n)$. The following
proposition is due to Cochran and Harvey (for example see
\cite[Proposition 6.1]{C}).

\begin{prop}
\label{prop:Seifert form} The following sequence is exact.
$$
H_1(\Sigma;\bbk_n)\ox_{\bbk_n}\bbk_n[t^{\pm 1}] \xrightarrow{d}
H_1(Y;\bbk_n)\ox_{\bbk_n}\bbk_n[t^{\pm 1}] \to \SA_n(K) \to 0
$$
where $d(\a\ox 1) = (i_+)_*\a\ox t - (i_-)_*\a\ox 1$.
\end{prop}

\begin{cor}\cite[Corollary 6.2]{C}
\label{cor:Seifert form} If the classical Alexander module of $K$
is not 1, then $\SA_n(K)$, $n>0$, has a square presentation matrix
of size $r = \text{max}\{0,-\chi(\Sigma)\}$ each entry of which is a
Laurent polynomial of degree at most 1. Specifically, we have the
presentation
$$
(\bbk_n[t^{\pm 1}])^r \xrightarrow{\partial} (\bbk_n[t^{\pm 1}])^r
\to \SA_n(K) \to 0
$$
where $\partial$ arises from the above proposition. If $n=0$, then
the same holds with $r$ replaced by $\b_1(\Sigma)$.
\end{cor}

\begin{defn}
\label{defn:Seifert form} The above presentation matrix is called
an  \emph{order $n$ localized Seifert presentation matrix} for
$K$.
\end{defn}

\begin{prop}\label{prop:failure} Theorem~\ref{thm:main} is false
for $n\geq 2$ if the hypothesis on the degree of the Alexander
polynomial is weakened to allow knots whose Alexander polynomials
have degree $2$.
\end{prop}
\begin{proof}
Let $K$ denote the knot shown in Figure~\ref{fig:ribbon}, well
known as the ``simplest'' ribbon knot.
\begin{figure}[htbp]
\setlength{\unitlength}{1pt}
\begin{picture}(191,156)
\put(0,0){\includegraphics{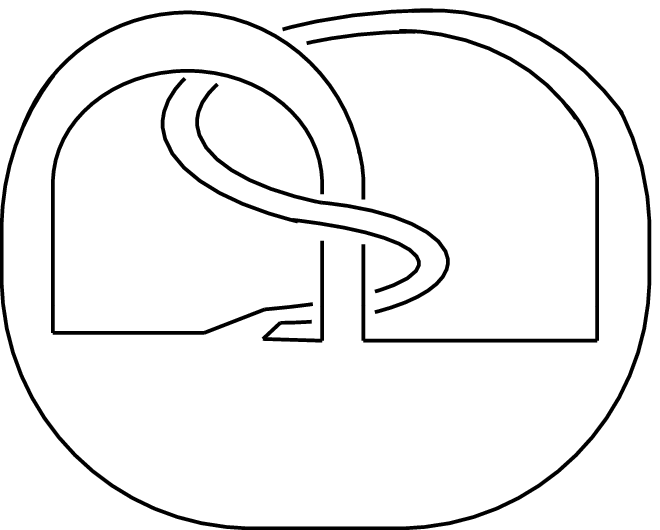}}
\end{picture}
\caption{}\label{fig:ribbon}
\end{figure}

We claim that the conclusions of Theorem~\ref{thm:main} for $n\geq
2$ are false for $K$ because, loosely speaking, any knot with the
same classical Alexander module and first higher-order Alexander
module as $K$ is topologically slice by recent work of Friedl and
Teichner. Details follow.

One easily checks from direct calculation that if $\Delta$ is one
of the two obvious ribbon disks obtained by ``cutting one of the
bands'' in the Figure, then $G=\pi_1(B^4-\Delta)$ is isomorphic to
$\mathbb{Z}[\frac{1}{2}]\rtimes \mathbb{Z}$ and in particular is
$1$-solvable, that is $G^{(2)}=0$. Let $P$ denote $\pi_1(M_K)$
where $M_K$ is the zero surgery on $K$. Then the inclusion map
induces an epimorphism $\phi:P\to G$. Since $G^{(3)}=0$, this
factors through $P/P^{(3)}$.

Suppose the conclusions of the main theorem were true for $K=K_0$
and let $K_1$ denote a knot other than $K$ that satisfies the
conclusions for some $n\geq 2$. By property $(1)$ of
Theorem~\ref{thm:main} $K_1$ is $(0)$-solvable and hence has Arf
invariant zero (this also follows from property $(4)$). Let $P_1$
denote $\pi_1(M_{K_1})$ where $M_{K_1}$ is the zero surgery on
$K_1$. As a consequence of property $(3)$, there is an isomorphism
$f: P_1/(P_1)^{(3)}\to P/P^{(3)}$ (since the conjugacy classes of
the longitudes correspond under the isomorphism given by $(3)$).
Hence there is an epimorphism $\phi_1:P_1\to G$ factoring through
$f$ (essentially $\phi \circ f$). Recall that for any coefficient
system $\psi:Q\to G$, where $Q$ is a group, we have the
isomorphism
$$
H_1(Q;\mathbb{Z}G)\cong \frac{\ker\psi}{[\ker\psi,\ker\psi]}
$$
and so in particular
$$
H_1(P_1;\mathbb{Z}G)\cong
\frac{\ker\phi_1}{[\ker\phi_1,\ker\phi_1]} \text{ and }
H_1(P;\mathbb{Z}G)\cong \frac{\ker\phi}{[\ker\phi,\ker\phi]}.
$$
Since $G^{(2)}=0$, $(P_1)^{(2)}\subset \ker\phi_1$ and so
$(P_1)^{(3)}\subset  [\ker\phi_1,\ker\phi_1]$. Similarly for $P$.
It follows that
$$
H_1(P_1;\mathbb{Z}G)\cong H_1(P_1/(P_1)^{(3)};\mathbb{Z}G)\cong
 H_1(P/(P)^{(3)};\mathbb{Z}G)\cong H_1(P;\mathbb{Z}G)
$$
Since $K_1$ and $K$ are nontrivial, $M_{K_1}$ and $M_K$ are
aspherical and so
$$ H_1(M_{K_1};\mathbb{Z}G)\cong
H_1(P_1;\mathbb{Z}G) \text{ and } H_1(M_{K};\mathbb{Z}G)\cong
H_1(P;\mathbb{Z}G).
$$
In summary, $H_1(M_{K_1};\mathbb{Z}G)\cong H_1(M_{K};\mathbb{Z}G)$
and thus it follows from \cite[Theorem 8.1]{FT} that $K_1$ is a
(topologically) slice knot, contradicting property $(2)$ of
Theorem~\ref{thm:main}.
\end{proof}

\begin{figure}[htbp]
\setlength{\unitlength}{1pt}
\begin{picture}(262,71)
\put(10,37){$\eta_1$} \put(120,37){$\eta_m$} \put(52,39){$\dots$}
\put(206,36){$\dots$} \put(183,37){$J_1$} \put(236,38){$J_m$}
\put(174,9){$K(\eta_1,\dots,\eta_m,J_1,\dots,J_m)$}
\put(29,7){$K$} \put(82,7){$K$}
\put(20,20){\includegraphics{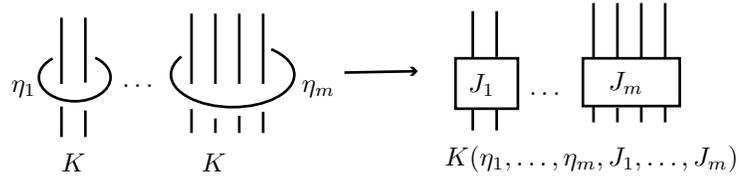}}
\end{picture}
\caption{$K(\eta_1,\dots,\eta_m,J_1,\dots,J_m)$: Genetic
modification of $K$ by $J_i$ along $\eta_i$}\label{fig:infection}
\end{figure}

We explain our method to construct desired examples. This
technique is called \emph{genetic modification} which results in a
satellite of a knot. One can find a detailed explanation of
genetic modification in \cite{COT2}. We briefly review this
construction. Let $K$ be a knot and $\{\eta_1,\eta_2,\ldots,
\eta_m\}$ be an oriented trivial link in $S^3$ which misses $K$.
Suppose $\{J_1,J_2,\ldots, J_m\}$ is an $m$-tuple of auxiliary
knots. For each $\eta_i$, remove a tubular neighborhood of
$\eta_i$ in $S^3$ and glue in a tubular neighborhood of $J_i$
along their common boundary, which is a torus, in such a way that
the longitude of $\eta_i$ is identified with the meridian of $J_i$
and the meridian of $\eta_i$ with the longitude of $J_i$. The
resulting knot is denoted by $K(\eta_i,J_i)$ and called the result
of \emph{genetic modification performed on the seed knot $K$ with
the infection $J_i$ along the axis $\eta_i$}. This construction
can also be described in the following way. For each $\eta_i$,
take an embedded disk in $S^3$ bounded by $\eta_i$ such that it
meets with $K$ transversally. Cut off $K$ along the disk, grab the
cut strands of $K$, tie them into the knot $J_i$ with 0-framing as
in Figure~\ref{fig:infection}.

Compare the following proposition with \cite[Proposition
3.1]{COT2}.

\begin{prop}
\label{prop:genetic modification} Let $K$ be a knot and $M$ be
zero surgery on $K$. Suppose $n\geq 1$. Suppose $[\eta_i]\in
\1(M)\2n$ and $J_i$ has vanishing Arf invariant for $i=1,2,\ldots,
m$. Then $K$ is $n$-solvequivalent to $K'\equiv K(\eta_i,J_i)$.
\end{prop}
\begin{proof}
Since the $J_i$ have vanishing Arf invariant, they are $0$-solvable
(see \cite[Remark 1.3.2]{COT1}). Let $W_i$ be a 0-solution for
$J_i$. By doing surgery along $\1(W_i)^{(1)}$  we may assume
$\1(W_i) \cong \bbz$, generated by the meridian of $J_i$.

Denote zero surgery on $K'$ in $S^3$ by $M'$. We construct an
$n$-cylinder between $M$ and $M'$. Take $M\times[0,1]$. Note that
$\partial W_i = M_i = (S^3\setminus N(J_i)) \cup_{S^1\times S^1}
S^1\times D^2$ where $N(J_i)$ denotes a tubular neighborhood of
$J_i$ in $S^3$ and $\{*\}\times D^2$ in $S^1\times D^2$ denotes
the surgery disk. Take the union of $M\times [0,1]$ and $W_i$'s in
such a way that we identify a product neighborhood of
$\eta_i\times \{1\}$ in $M\times \{1\}$, with the $i^{th}$ copy of
$S^1\times D^2$ in $\partial W_i$ for each $i$. We denote the
4-manifold resulted from this construction by $W$. One easily sees
that $\partial W = M \coprod -M'$. We claim that $W$ is an
$n$-cylinder. Using a Mayer-Vietoris sequence  ($M\times [0,1]$
union $W_i$ intersecting along $\eta_i\times D^2$), one easily
sees that $H_1(W) \cong \bbz$ and the inclusions from $M$ to $W$
and $M'$ to $W$ induce isomorphisms on the first homology. Also
from the Mayer-Vietoris sequence one observes that $H_2(W) \cong
H_2(M)\op (\op^m_{i=1} H_2(W_i)) \cong H_2(M')\op (\op^m_{i=1}
H_2(W_i))$. Note that the generator of $H_2(M)$ under the map
$H_2(M) \lra H_2(W)$ maps to the image of the generator of
$H_2(M')$ under the map $H_2(M') \lra H_2(W)$ since they are
represented by capped-off Seifert surfaces. Thus
$H_2(W)/i_*(H_2(\partial W)) \cong \op^m_{i=1}H_2(W_i)$ where
$i_*$ is the homomorphism induced from the inclusion $i:
\partial W \lra W$. Recall that $\1(W_i)$ $(\cong \bbz)$ is
generated by the meridian of $J_i$ which is identified with the
longitude of $\eta_i$. Hence $\1(W_i)$ maps into $\1(W)^{(n)}$
since $[\eta_i]$ lies in $\1(M)^{(n)}$ by hypothesis and
$\1(M)^{(n)}$ maps into $\1(W)^{(n)}$. This implies that
0-surfaces in $W_i$ are $n$-surfaces in $W$. Now by naturality of
intersection forms, one can see that the union of the
0-Lagrangians (with 0-duals for $W_i$) ($0\le i \le m$)
constitutes an $n$-Lagrangian (with $n$-duals) for $W$. One checks
easily that $W$ is spin since each $W_i$ was spin.
\end{proof}

Suppose $K,\eta_i,J_i$ and $K'$ are as in
Proposition~\ref{prop:genetic modification}. Denote zero surgery
on $J_i$ in $S^3$ by $M_i$. Suppose $W$ is the $n$-cylinder with
$\partial W= M \coprod -M'$ constructed above. Suppose
$\psi:\1(W)\lra \G$ is a map to an arbitrary $n$-solvable PTFA
group $\G$. Let $\phi, \phi'$, and $\phi_i$ denote the induced
maps on $\1(M), \1(M')$, and $\1(M_i)$ respectively. For each
$J_i$, we denote by $\rho_\bbz(J_i)$ the $\rho$-invariant
$\rho(M_i,\zeta_i)$ where $\zeta_i : \1(M_i) \lra \bbz$ is the
abelianization. The following lemma reveals the additivity of the
$\rho$-invariant in our current setting.

\begin{lem}
\label{lem:additivity}
$$
\rho(M,\phi) - \rho(M',\phi') = \sum^m_{i=1}\e_i\rho_{\mathbb{Z}}(J_i)
$$
where $\e_i=0$ or $1$ according as $\phi(\eta_i)=e$ or not.
\end{lem}
\begin{proof}
It follows from \cite[Proposition 3.2]{COT2} that
$$
\rho(M,\phi) - \rho(M',\phi') = \sum^m_{i=1}\rho(M_i,\phi_i)
$$
Since $\pi_1(W_i)\cong \mathbb{Z}$, $\phi_i$ factors through
$\mathbb{Z}$ generated by $\eta_i$. Thus its image is zero if
$\phi(\eta_i)=e$, and $\mathbb{Z}$ if not ( recall a PTFA group is
torsion-free). In the former case, $\rho(M_i,\phi_i)=0$ by
property ($3$) of Proposition~\ref{prop:rho invariants}. In the
latter case, $\rho(M_i,\phi_i)=\rho_\bbz(J_i)$ by property ($2$)
of Proposition~\ref{prop:rho invariants}.
\end{proof}

The following two theorems are the key theorems that show the
existence of a knot that is $n$-solvequivalent to a given knot $K$
but not $(n.5)$-solvequivalent to $K$. The first significantly
generalizes \cite[Theorem 4.3]{CT} in two ways. Firstly, their
theorem applies only when $K$ is a genus $2$ fibered knot.
Secondly their theorem only covers the case when $\partial W=M$.
For a group $G$, let $G\2k_r$ be the \emph{$k$-th rational derived
group of $G$} (see \cite[Section 3]{Ha}).

\begin{thm}
\label{thm:injectivity} Let $n\in\bbn$. Let $K$ be a knot for
which the degree of the Alexander polynomial is greater than $2$
(if $n=1$, degree equal to $2$ is allowed). Let $M$ be the zero
framed surgery on $K$ in $S^3$. Suppose $\Sigma$ is a Seifert
surface for $K$. Then there exists an oriented trivial link
$\{\eta_1,\eta_2,\dots,\eta_m\}$ in $S^3-\Sigma$ that satisfies
the following:
\begin{itemize}
\item [(1)] $\eta_i\in\1(M)\2n$ for all $i$. Moreover, the
$\eta_i$ bound (smoothly embedded) symmetric capped gropes of
height $n$, disjointly embedded in $S^3-K$ (except for the caps,
which will hit $K$).

\item [(2)] For every $n$-cylinder $W$ with $M$ as one of its boundary
components, there exists {\bf some} $i$ such that
$j_*(\eta_i)\notin\1(W)\nn_{r}$ where $j_*:\1(M)\lra\1(W)$ is
induced from the inclusion $j:M \lra W$. The number of such $i$'s
is at least $\frac12 (d-2)$ if $n>1$ or at least $\frac12 d$ if
$n=1$ where $d$ denotes the degree of the Alexander polynomial of
$K$.
\end{itemize}
\end{thm}
\begin{proof}
Let $S \equiv \1(M)^{(1)}$. Suppose $W$ is an $n$-cylinder with
$M$ as one of its boundary components. The inclusion $j:M\lra W$
induces a map $j:S\lra\1(W)\21$. Let $G \equiv
\1(W)\21=\1(W)\21_{r}$. The last equality holds since $H_1(W)
\cong \bbz$. For convenience, let us use the same notation
$\Sigma$ for the capped-off Seifert surface. Let $g \equiv \frac12
\rk_\bbq H_1(M\setminus \Sigma;\bbq)$ (which is equal to the genus
of $\Sigma$). Let $F$ be the free group of rank $2g$. Choose a map
$F\lra \1(M\setminus \Sigma)$ that induces an isomorphism on
$H_1(F;\bbq)$. Let $i$ be the composition $F\lra \1(M\setminus
\Sigma)\subset S$. For a group $G$, let $G_k$, for $k\in \bbn_0$,
denote $G/G\2k_r$. By \cite[Corollary 3.6]{Ha} $G_k$ is PTFA. By
\cite[Proposition 2.5]{COT1} $\bbz G_k$ embeds into the (skew)
quotient field of fractions, which is denoted by $\bbk(G_k)$. By
\cite[Proposition II.3.5]{Ste} $\bbk(G_k)$ is flat over $G_k$.
Since $\1(W)/\1(W)\21 \cong \bbz$, if we let $\G_k \equiv
\1(W)/\1(W)_r^{(k+1)}$, there is a short exact sequence $1\to G_k
\to \G_k \xrightarrow{\pi} \bbz \to 1$ where $\pi$ is the
abelianization. Then we have a PID $\bbk(G_k)[t^{\pm 1}]$ such
that $\bbz\G_k \subset \bbk(G_k)[t^{\pm 1}] \subset \SK$ where
$\SK$ is the (skew) quotient field of fractions of $\bbz\G_k$
($\G_k$ is PTFA by \cite[Corollary 3.6]{Ha}, hence $\bbz\G_k$
embeds into the skew quotient field of fractions).

Apply Theorem~\ref{thm:paininbutt} to find a finite collection
$\SP_{n-1}$ of $2g-1$-tuples ($2g$ if $n=1$) of elements of
$F^{(n-1)}$. Since $S^{(n-1)}=\1(M)\2n$, the image of the union of
the elements of $\SP_{n-1}$ under $i : F \lra S$ is a finite set
$\{\a_1,\dots,\a_m\}$ of elements of $\1(M)\2n$ as required by the
part (1) in the statement. By Proposition~\ref{prop:geomalg} the
induced map $j:S\lra G$ is an algebraic $n$-solution, hence an
algebraic $(n-1)$-solution (see Remark~\ref{rem:algebraic}(3)). By
our choice of $\SP_{n-1}$ from Theorem~\ref{thm:paininbutt}, at
least one tuple $\{w_1,w_2,\ldots ,w_{2g-1}\}\in\SP_{n-1}$
($\{w_1,w_2,\ldots ,w_{2g}\}$ if $n=1$) maps to a generating set
of $H_1(S;\bbk(G_{n-1}))$. Notice that $H_1(S;\bbk(G_{n-1})) \cong
H_1(M;\bbk(G_{n-1})[t^{\pm 1}])$ and $H_1(G;\bbk(G_{n-1})) \cong
H_1(W;\bbk(G_{n-1})[t^{\pm 1}])$ as $\bbk(G_{n-1})$-modules. By
Theorem~\ref{thm:rank} (with $\G = \1(W)/\1(W)\2n_r, [\G,\G] =
\1(W)\21/\1(W)\2n_r = G_{n-1}$), we see that at least $\frac12
(d-2)$ of $\{w_1,...,w_{2g-1}\}$ (if $n>1$) or at least $\frac12
d$ of $\{w_1,...,w_{2g}\}$ (if $n=1$) map non-trivially under
$F^{(n-1)}\to S^{(n-1)}\overset{j}{\to}G^{(n-1)}_{r}/G\2n_{r}$
$(\cong H_1(G;\bbz G_{n-1})$ modulo $\mathbb{Z}$-torsion). Hence
at least $\frac12 (d-2)$ (if $n>1$) or at least $\frac12 d$ (if
$n=1$) of the $\a_i$ have the property that $j_*(\a_i)\notin
G\2n_{r}=(\1(W))^{\nn}_{r}$. Since each $\a_i$ actually lies in
$\1(S^3\setminus \Sigma)^{(n)}$ we can represent the $\a_i$ by
simple closed curves in the complement of the chosen Seifert
surface $\Sigma$ for the knot $K$. We can alter these by crossing
changes until the collection of $\a_i$ forms a trivial link in
$S^3$. Moreover, by \cite[Lemma 3.8]{CT} and its proof, we may
choose representatives in the same homotopy classes that have the
much stronger property that they bound (smoothly embedded)
symmetric capped gropes of height $n$, disjointly embedded in
$S^3-K$ except for the caps. This is the collection $\{\eta_i\}$
required.

\end{proof}

Let $c_M$ denote the universal bound for all $\rho$-invariants of
a fixed 3-manifold $M$ given by \cite[Theorem 4.10]{ChG} and
Ramachandran \cite[Theorem 3.1.1]{R}. That is, $|\rho(M,\phi)| <
c_M$ for every representation $\phi : \1(M) \lra \G$ where $\G$ is
a group. For a detailed discussion see \cite{CT}. The final
conclusion of the following theorem was proved in \cite{CT} in the
case that $K$ is itself an $n$-solvable knot.
\begin{thm}
\label{thm:infection} Let $n\in \bbn$. Let $K$ be a knot in $S^3$.
Suppose $\{\eta_1,\eta_2,\dots,\eta_m\}$ is an oriented trivial
link in $S^3$ that misses $K$ and has properties ($1$) and ($2$)
of Theorem~\ref{thm:injectivity}. Then, for every $m$-tuple
$\{J_1,J_2,\dots,J_m\}$ of Arf invariant zero knots for which
$\rho_\bbz(J_i)> 2c_M$, the knot $K' \equiv K(\eta_i,J_i)$ (formed
by genetic modification) is $n$-solvequivalent to $K$ but not
$(n.5)$-solvequivalent to $K$. Moreover $K'-K$ is of infinite
order in $\ff$. (In particular, $K'\#(-K)$ is $(n)$-solvable.)
\end{thm}
\begin{proof}
By Proposition~\ref{prop:genetic modification}, $K'$ is
$n$-solvequivalent to $K$. Let $M$ be zero surgery on $K$ in $S^3$
and $M'$ zero surgery on $K'$ in $S^3$. Let $M_i$ be zero surgery
on $J_i$ in $S^3$, $1\le i\le m$. Let $W$ be an $n$-cylinder with
$\partial W = M \coprod -M'$ constructed as in the proof of
Proposition~\ref{prop:genetic modification}. Suppose $M'$ is
$(n.5)$-solvequivalent to $M$ via $V$. We will show that this
leads us to a contradiction. Let $X\equiv W\cup_{M'} V$. Thus
$\partial X = M\coprod -M$. We assume $\partial_+X = M$ and
$\partial_-X = -M$. Since $V$ is an $(n.5)$-cylinder, it is an
$n$-cylinder. As in Proposition~\ref{prop:transitivity}, one can
show that $X$ is an $n$-cylinder. Let $\G \equiv
\1(X)/\1(X)\nn_r$, an $n$-solvable PTFA group. Let $\psi : \1(X)
\lra \G$ be the projection. Let $\phi_+$, $\phi_-$, and $\phi'$
denote the induced maps on $\1(\partial_+X)$, $\1(\partial_-X)$,
and $\1(M')$, respectively. By Lemma~\ref{lem:additivity} we have
$$
\rho(M,\phi_+)-\rho(M',\phi')=\sum^m_{i=1}\e_i\rho_{\bbz}(J_i)
$$
where $\e_i=0$ or $1$ according as $\phi_+(\eta_i)=e$ or not. On
the other hand, since $V$ is an $(n.5)$-cylinder,
Theorem~\ref{thm:rho=0} applies to say that
$$\rho(M',\phi')-\rho(M,\phi_-)=0.
$$
Hence
$$
\rho(M,\phi_+)-\rho(M,\phi_-)=\sum^m_{i=1}\e_i\rho_\bbz(J_i).
$$
Note that since $X$ is an $n$-cylinder and the collection
$\{\eta_i\}$ was chosen to satisfy property $(2)$ of
Theorem~\ref{thm:injectivity}, there exists at least one $i$ such
that $\psi(\eta_i)\neq e$ from which it follows that
$\phi_+(\eta_i)\neq e$. Thus
$$
\rho(M,\phi_+) - \rho(M,\phi_-) \ge \rho_{\bbz}(J_i)> 2c_M.
$$
which is a contradiction. Therefore $K'$ is not
$(n.5)$-solvequivalent to $K$.

We will now show that $K'-K\in \SF_{(n)}$, that is, $K'\#(-K)$ is
 $n$-solvable. Note that since $K'$ was obtained from $K$ by
genetic modification along the $\eta_i$, $K'\#(-K)$ can be
obtained from $K\#(-K)$ by performing a genetic modification along
the ``same'' axes $\eta_i$. It is only necessary to observe that
we can perform the connected sum in such a way as to preserve the
fact that the $\eta_i$ lie in the $n^{th}$-derived group of
$\pi_1$. This is clear because $S^3-K'$ can be viewed as a
\emph{subspace} of $S^3-(K\#(-K))$. Notice that $K\#(-K)$ is a
slice knot, hence $n$-solvable. By \cite[Proposition 3.1]{COT2},
performing a genetic modification on an $n$-solvable knot using
Arf invariant zero knots $J_i$, along axes that lie in the
$n^{th}$-derived group results in another $n$-solvable knot. Hence
$K'\#(-K)$ is $n$-solvable.

Suppose $K'-K$ were of order $k > 0$ in $\ff$. We will show that
this yields a contradiction, implying that $K'-K$ is of infinite
order. Our assumption is equivalent to the fact that
$(\#^k_{i=1}K')\#(\#^k_{j=1}-K)$ is $(n.5)$-solvable. Let $N$
denote the zero surgery on $(\#^k_{i=1}K')\#(\#^k_{j=1}-K)$. Let
$V$ be an $(n.5)$-solution for $N$ and let $W$ be as above. Let
$M'_i$ be the $i$-th copy of $M'$ and $M_j$ be the $j$-th copy of
$M$, $1\le i,j\le k$. Take a standard (spin) cobordism $C$ from
$(\coprod^k_{i=1}M'_i)\coprod (\coprod^k_{j=1}-M_j)$ to $N$ which
is obtained from $((\coprod^k_{i=1}M'_i)\coprod
(\coprod^k_{j=1}-M_j)) \times [0,1]$ by adding $(2k-1)$ 1-handles
and then $(2k-1)$ 2-handles (see \cite[Lemma 4.2]{COT2} for more
detail). Let $X \equiv (\coprod^k_{i=1}W_i)\cup_{\coprod M'_i}C
\cup_{N} V$ where $W_i$ is the $i$-th copy of $W$. Then $\partial
X = (\coprod^k_{i=1}M_i) \coprod (\coprod^k_{j=1}-M_j)$. Let $\G
\equiv \1(X)/\1(X)\nn_r$, an $n$-solvable PTFA group. Let $\psi
:\1(X) \lra \G$ be the projection. Let $\phi_j$, $\phi'_i$,
$\phi^+_i$, and $\phi$ denote the restrictions of $\psi$ to
$\1(-M_j)$, $\1(M'_i)$, $\1(M_i)$ (the upper boundary of $W_i$),
and $\1(N)$ respectively. One sees that $C$ is an $(n.5)$-cylinder
since $H_2(C)/i_*(H_2(\partial C)) =0$. Thus by
Theorem~\ref{thm:rho=0}
$$
\sum^k_{i=1}\rho(M'_i,\phi'_i) + \sum^k_{j=1}\rho(-M_j,\phi_j) -
\rho(N,\phi)= 0.
$$
But since $V$ is an $(n.5)$-solution for $N$, $\rho(N,\phi) = 0$
by Theorem~\ref{thm:rho=0}. Therefore
$\sum^k_{i=1}\rho(M'_i,\phi'_i) = \sum^k_{j=1}\rho(M_j,\phi_j)$.
By Lemma~\ref{lem:additivity} (applied for each $i$ and then
summed),
$$
\sum^k_{i=1}\rho(M_i,\phi^+_i)  - \sum^k_{i=1}\rho(M'_i,\phi'_i) =
\sum^k_{i=1}\sum^m_{l=1}\e_{il}\rho_\bbz(J_l)
$$
where $\e_{il} = 0$ if $\phi^+_i(\eta_{il}) = e$ and $\e_{il} = 1$
if $\phi^+_i(\eta_{il}) \ne e$. (By $\eta_{il}$ we mean the $l$-th
axis $\eta_{l}$ for the $i$-th copy of $K$.) Thus
$$
\sum^k_{i=1}\rho(M_i,\phi^+_i)  - \sum^k_{j=1}\rho(M_j,\phi_j) =
\sum^k_{i=1}\sum^m_{l=1}\e_{il}\rho_\bbz(J_l).
$$
Now we claim that $X$ is an $n$-cylinder. Using a Mayer-Vietoris
sequence and \cite[Lemma 4.2]{COT2}, one can show that
$H_2(X)/i_*(H_2(\partial X)) \cong (\op^k_{i=1}H_2(W_i)) \op
H_2(V)$ and the union of the  $0$-Lagrangians with 0-duals for
$W_i$ ($1\le i\le k$) and the $n$-Lagrangian with $n$-duals for
$V$ constitutes an $n$-Lagrangian with $n$-duals for $X$. Since,
for any $i$, $X$ is an $n$-cylinder with $M_i$ as one of its
boundary components, by Theorem~\ref{thm:injectivity}, for each
$i$ there exists some $l_i$ such that $\psi(\eta_{il_i}) \ne e$.
Hence $\phi^+_i(\eta_{il_i})\ne e$. Then the right-hand side of
the last equation is greater than $2kc_M$, while the left-hand
side is less than $2kc_M$, a contradiction.
\end{proof}

Finally we prove the main theorem.
\begin{proof}[{\bf Proof of Theorem~\ref{thm:main} (Main Theorem)}]
Let $K$ be a knot for which the degree of the Alexander polynomial
is greater than 2 (if $n=1$, degree equal to $2$ is allowed). Let
$M$ denote the zero surgery on $K$. Apply
Theorem~\ref{thm:injectivity} to choose an oriented trivial link
$\{\eta_1,\eta_2,\ldots, \eta_m\}$ in $S^3$ that misses $K$
(except for the caps, which will hit $K$), and such that the
$\eta_i$ bound height $n$ (smoothly embedded) symmetric capped
gropes disjointly embedded in $S^3-K$. Let $c_M$ be the universal
bound for $\rho$-invariants for $M$ as mentioned above. Choose
knots $J^i$ ($i\in \bbn)$ with vanishing Arf invariants
inductively as follows. First, choose $J^1$ such that
$\rho_\bbz(J^1) > 2c_M$. Suppose $J^{k-1}$ has been constructed.
Then choose $J^k$ such that $\rho_\bbz(J^k) > 2c_M +
2m\rho_\bbz(J^{k-1})$. These $J^i$ are easily found by taking the
connected sum of a suitably large even number of copies of the
left-handed trefoil. (Note that for the left-handed trefoil $J$,
$\rho_\bbz(J)=4/3$ by Proposition~\ref{prop:rho invariants}(4).)
However, to achieve the final conclusion of Part $(1)$, we must
choose each $J^i$ to be a suitably large connected sum of the knot
shown in \cite[Figure 1.7]{CT}. This knot has the same $\rho_\bbz$
as the left-handed trefoil knot \cite[Lemma 4.4]{CT}. For each
$i$, let $J^i_j$ be the $j$-th copy of $J^i$, $1\le j \le m$. Now
define $K_i \equiv K(\eta_1,\dots,\eta_m,J^i_1,\dots,J^i_m)$, the
result of genetic modification performed on $K$ with the
infections $J^i_j$ along the axes $\eta_j$ $(1\le j\le m)$. Set
$K_0 \equiv K$.

\noindent {\bf Part (1)} : It follows from the last sentence of
Theorem~\ref{thm:infection} that $K_i-K \in \SF_{(n)}$, $i>0$ and
hence that $K_i$ is $n$-solvequivalent to $K$ (or apply
Proposition~\ref{prop:genetic modification}). To show that $K_i$
and $K$ cobound, in $S^3\times\[0,1\]$, a (smoothly embedded)
embedded symmetric Grope of height $n+2$, merely apply the proof
of \cite[Theorem 3.7]{CT}.

\noindent {\bf Part (2)} : Suppose $i>j\ge 0$. Let $M_i$ be the
zero surgery on $K_i$ in $S^3$. Suppose $M_i$ is
$(n.5)$-solvequivalent to $M_j$ via $U$ ($\partial U = M_i\coprod
-M_j$). Suppose $V$ be an $n$-cylinder with $\partial V = M
\coprod -M_i$ and $W$ an $n$-cylinder with $\partial W =
M_j\coprod -M$ such that $V$ and $W$ are constructed as in the
proof of Proposition~\ref{prop:genetic modification}. Let $X
\equiv V\cup_{M_i}U \cup_{M_j} W$. Then $X$ is an $n$-cylinder
with $\partial_+ X = M$ and $\partial_-X = -M$ (see
Proposition~\ref{prop:transitivity}). Let $\G \equiv
\1(X)/\1(X)\nn_r$, an $n$-solvable PTFA group. Let $\psi : \1(X)
\lra \G$ be the projection. Let $\phi_+$, $\phi_i$, $\phi_j$, and
$\phi_-$ denote the restrictions of $\psi$ to $\1(\partial_+ X) (
= \1(M))$, $\1(M_i)$, $\1(M_j)$, and $\1(\partial_-X) ( = \1(-M))$
respectively. Since $U$ is an $(n.5)$-cylinder, by
Corollary~\ref{cor:rho=0}, $\rho(M_i,\phi_i) = \rho(M_j,\phi_j)$.
On the other hand, by Lemma~\ref{lem:additivity},
$$
\rho(M,\phi_+) - \rho(M_i,\phi_i) =
\sum^m_{k=1}\e_k\rho_\bbz(J^i_k)
$$
where $\e_k = 0$ if $\phi_+(\eta_k) = e$ and $\e_k = 1$ if
$\phi_+(\eta_k) \ne e$. Similarly,
$$
\rho(M,\phi_-) - \rho(M_j,\phi_j) =
\sum^m_{k=1}\e'_k\rho_\bbz(J^j_k)
$$
where $\e'_k = 0$ if $\phi_-(\eta_k) = e$ and $\e'_k = 1$ if
$\phi_-(\eta_k) \ne e$. Thus we have
$$
\rho(M,\phi_+) - \rho(M,\phi_-) =
\sum^m_{k=1}\e_k\rho_{\bbz}(J^i_k) -
\sum^m_{k=1}\e'_k\rho_{\bbz}(J^j_k).
$$
Since $J^j_k$ is merely a copy of $J^j$,
$\rho_\bbz(J^j_k)=\rho_\bbz(J^j)$, and similarly
$\rho_{\bbz}(J^i_k)=\rho_\bbz(J^i)$. Moreover, by
Theorem~\ref{thm:injectivity} applied to $M=\partial_+ X$, there
is some $k$ such that $\phi_+(\eta_k) = \psi(\eta_k) \ne e$. Thus
the right hand side above is greater than
$\rho_\bbz(J^i)-m\rho_\bbz(J^j)$, which, by our choice of $J^i$
and $J^j$ is greater than $2c_M$. This is a contradiction since
the left-hand side above has absolute value less than $2c_M$. To
show that $K_i$ and $K_j$ do not cobound, in $S^3\times\[0,1\]$,
an embedded symmetric Grope of height $(n+2.5)$, merely note that
the proof of \cite[Theorem 8.11]{COT1} clearly applies to show
that if they did bound such a Grope then they would be
$(n.5)$-solvequivalent.

\noindent {\bf Part (3)} : This follows from \cite[Theorem
8.1]{C}.

\noindent {\bf Part (4)} : We shall employ the terminology set out
just above Proposition~\ref{prop:Seifert form}. Let $K_*$ denote
one of the $K_i$ and let $G_*$ denote its knot group. Let $G$
denote the knot group of $K$. By part $(3)$, these two knots share
the same higher-order (integral) Alexander modules up to order
$n-1$ and $G/G^{(i+1)}$ is isomorphic to $G_*/(G_*)^{(i+1)}$ for
all $i\leq n$. Therefore in considering the $i^{th}$ higher-order
Alexander modules of these knots, we may consider that they are
modules over the same ring, $\mathbb{Z}\G_i\equiv
\mathbb{Z}[G/G^{(i+1)}]$, as long as $i\leq n$. We shall show that
there exist Seifert surfaces $\Sigma_*$ for $K_*$ and $\Sigma$ for
$K$ with respect to which the $i^{th}$ order Seifert
presentations, $i_{\pm}:H_1(\Sigma;\mathbb{Z}\G'_i)\to
H_1(Y;\mathbb{Z}\G'_i)$ and
$i_{\pm}:H_1(\Sigma_*;\mathbb{Z}\G'_i)\to
H_1(Y_*;\mathbb{Z}\G'_i)$ are identical up to isomorphisms
identifying the domain and range of each, as long as $i\leq n-1$.
(Recall that $\G_i'$ is the commutator subgroup of $\G_i$.) A
Seifert presentation determines the map $d$ in
Proposition~\ref{prop:Seifert form} and hence there are bases with
respect to which the localized $i^\textrm{th}$-order Seifert
matrices for $K_*$ and $K$ are identical. We may assume that
$n\geq1$.

Let $E(K)$ denote the exterior of $K$ in $S^3$. The continuous map
$f:E(K_*)\to E(K)$ that induces all of these isomorphisms is
described as follows (see \cite[Theorem 8.1]{C}). Recall that
$E(K_*)$ is constructed from $E(K)$ by replacing a collection of
solid tori $\eta_j \times D^2$ by a collection of knot exteriors
$E(J_j)$. Since it is well known that there is always a degree one
map, $f_j$, relative boundary  from $E(J_j)$ to $E(unknot)\equiv
\eta_j \times D^2$, there is a degree one map relative boundary,
$f$, from $E(K_*)$ to $E(K)$. Recall also that $\eta_j\in
\pi(M)^{(n)}$ by choice. But in fact, the $\eta_j$ actually
produced by Theorem~\ref{thm:paininbutt}, lie in $F^{(n-1)}$ where
$F\to \pi_1(M-\Sigma^*)$ where $\Sigma^*$ is a capped-off Seifert
surface. Hence, the circles $\eta_j$ can be chosen to miss the
given Seifert surface $\Sigma$ and actually represent elements of
$\pi_1(Y)^{(n-1)}\subset G^{(n)}$. Therefore we can use the `same'
Seifert surface for $K_*$ as for $K$. For simplicity for the rest
of this proof we shall assume that there is just one circle,
$\eta$. The proof is no different for a collection. Note that $f$
is the `identity' on a neighborhood of $\Sigma_*$ (mapping to
$\Sigma$) and restricts to a degree one map relative boundary
$f:Y_*\to Y$. We need only show that this map carries
$H_1(Y_*;\mathbb{Z}\G'_i)$ isomorphically to
$H_1(Y;\mathbb{Z}\G'_i)$. Briefly, this is true because $\eta\in
G^{(n)}$ and thus goes to zero in $\G_i$ if $i\leq n-1$. It
follows that the entire group $\pi_1(E(J))$, where $J$ is the
infection knot, maps to zero in $\G'_i$. Thus in a Mayer-Vietoris
analysis, $H_1(E(J);\mathbb{Z}\G'_i)$ really has untwisted
coefficients and thus is not distinguishable from the case that
$J$ is a trivial knot. But if $J$ were trivial then $K_*=K$ and
$Y_*=Y$. For more details, the proof is the same as the proof of
\cite[Theorem 8.2]{C} where it is shown that each of
$H_1(Y_*;\mathbb{Z}\G_i)$ and $H_1(Y;\mathbb{Z}\G_i)$ is
isomorphic to the quotient of $H_1(Y\setminus (\eta \times
D^2);\mathbb{Z}\G_i)$ by the submodule generated by the meridian
of $\eta$.

\noindent {\bf Part (5)} : For $j=0$ this follows directly from
Theorem~\ref{thm:infection}. For general $j$ the proof is
essentially the same as the proof of part $(2)$. We outline the
proof. Here $i$ and $j$ are fixed.

Let $M_i$, $M_j$, and $M$
denote the zero surgeries on $K_i$, $K_j$, and $K$ as usual.
Suppose $K_i-K_j$ were of order $k>0$, that is to say
$(\#^k_{l=1}K_i)\#(\#^k_{l=1}-K_j)$ is $(n.5)$-solvable. Then by
Proposition~\ref{prop:old and new}, $M_{\#K_i}$ is
$(n.5)$-solvequivalent to $M_{\#K_j}$ via some $U$ as in the proof
of part $(2)$. Refer to the schematic Figure~\ref{fig:cobordism} below.
\begin{figure}[h]
\setlength{\unitlength}{1pt}
\begin{picture}(118,203)
\put(58,87){$U$}
\put(58,116){$C$}
\put(58,57){$D$}
\put(14,155){$V^1$}
\put(58,155){$V^l$}
\put(105,155){$V^k$}
\put(13,20){$W^1$}
\put(58,20){$W^l$}
\put(102,20){$W^k$}
\put(108,71){$-M_{\#K_j}$}
\put(110,106){$M_{\#K_i}$}
\put(-6,143){$M^1_i$}
\put(81,143){$M_i^l$}
\put(35,170){$\dots$}
\put(35,11){$\dots$}
\put(79,11){$\dots$}
\put(79,170){$\dots$}
\put(14,177){$M^{1_v}$}
\put(57,177){$M^{l_v}$}
\put(100,177){$M^{k_v}$}
\put(123,143){$M_i^k$}
\put(76,35){$-M_j^l$}
\put(123,35){$-M_j^k$}
\put(-14,35){$-M_j^1$}
\put(10,-2){$-M^{1_w}$}
\put(53,-2){$-M^{l_w}$}
\put(98,-2){$-M^{k_w}$}
\put(10,10){\includegraphics{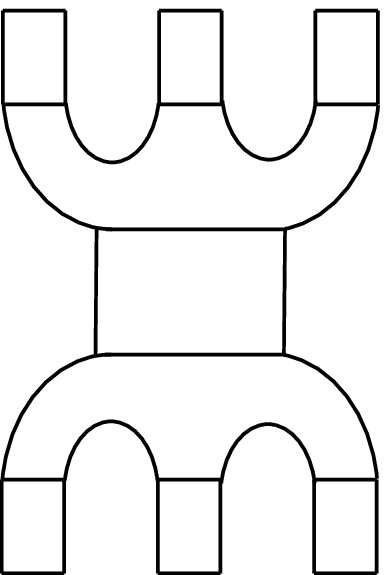}}
\end{picture}
\caption{}\label{fig:cobordism}
\end{figure}

Now let $C$ be the standard cobordism between
$M_{\#K_i}$ and the disjoint union of $k$ copies of $M_i$, which
we denote by $M_i^l$ for $1\leq l\leq k$, and let $D$ be the
standard cobordism between $-M_{\#K_j}$ and the disjoint union of
$k$ copies of $-M_j$, which we denote by $-M_j^l$ for $1\leq l\leq
k$. Let $V^l$ denote the $l^{th}$ copy of the standard n-cylinder
between $M_i$ and $M$, constructed in the proof of
Proposition~\ref{prop:genetic modification}, so that $\partial
V^l=M^{l_v}\coprod -M_i^l$ where we use $M^{l_v}$ to denote the
copy of $M$ that occurs in the boundary of $V^l$. Similarly let
$W^l$ denote the $l^{th}$ copy of the standard $n$-cylinder
between $-M_j$ and $-M$ so that $\partial W^l=M^{l_w}\coprod
-M_j^l$ where we use $M^{l_w}$ to denote the copy of $M$ that
occurs in the boundary of $W^l$. Following the proof of Part
$(2)$, let
$$
X\equiv (\coprod_{l=1}^k V^l)\cup C\cup_{M_{\#K_i}}U
\cup_{M_{\#K_j}}D\cup (\coprod_{l=1}^k W^l).
$$
Recall that we saw in the proof of Theorem~\ref{thm:injectivity}
that $C$ and $D$ are $(n.5)$-cylinders since
$H_2(C)/i_*(H_2(\partial C))=0$ and $H_2(D)/i_*(H_2(\partial
D))=0$. It follows that $X$ is an $n$-cylinder with $\partial_+ X
= \coprod_{l=1}^k M^{l_v}$ and $\partial_-X = \coprod_{l=1}^k
-M^{l_w}$ (see also Proposition~\ref{prop:transitivity}). Let $\G
\equiv \1(X)/\1(X)\nn_r$, an $n$-solvable PTFA group. Let $\psi :
\1(X) \lra \G$ be the projection. Let $\phi^{l_v}$ and
$\phi^{l_w}$ denote the restrictions of $\psi$ to $ \1(M^{l_v})$,
and $\1(M^{l_w})$ respectively. Since $U$ is an $(n.5)$-cylinder
and otherwise $H_2(X)$ comes from the $V^l$ and the $W^l$, we
arrive at an expression similar to that of the proof of Part
$(2)$, except summed from $l=1$ to $k$:
$$
\sum_{l=1}^k(\rho(M^{l_v},\phi^{l_v}) - \rho(M^{l_w},\phi^{l_w}))
=\sum_{l=1}^k( \sum^m_{s=1}\e_{sl}\rho_{\bbz}(J^i) -
\sum^m_{s=1}\e'_{sl}\rho_{\bbz}(J^j)),
$$
where $\e_{sl} = 0$ if $\phi^{l_v}(\eta_{sl}) = e$ (where
$\eta_{sl}$ means the $s$-th axis $\eta_s$ for $K$ for
constructing the $l$-th copy of $K_i$) and $\e_{sl} = 1$ if
$\phi^{l_v}(\eta_{sl}) \ne e$, and where $\e'_{sl} = 0$ if
$\phi^{l_w}(\eta'_{sl}) = e$ (where $\eta'_{sl}$ means the $s$-th
axis $\eta_s$ for $K$ for constructing the $l$-th copy of $K_j$)
and $\e'_{sl} = 1$ if $\phi^{l_w}(\eta'_{sl}) \ne e$. By property
$(2)$ of Theorem~\ref{thm:injectivity} applied, for \emph{each}
$l$, to $M^{l_v}$, there is some $s$ such that
$\phi^{l_v}(\eta_{sl}) = \psi(\eta_{sl}) \ne e$. Thus, for each
$l$, the $l^{th}$ term of the summation on the right hand side
above is greater than $\rho_\bbz(J^i)-m\rho_\bbz(J^j)$, which, by
our choice of $J^i$ and $J^j$ is greater than $2c_M$. Thus the
right hand side is greater than $2kc_M$. This is a contradiction
since the left-hand side above has absolute value less than
$2kc_M$, since $M^{l_v}\cong M^{l_w}\cong M$.

\noindent {\bf Part (6)} : Let $g$ be the genus of $K$. Note that
in the construction of $K_i$ we could have chosen \emph{any}
Seifert surface for $K$, in particular one of genus $g$. Recall
that the circles $\eta_i$ are chosen in the complement of this
Seifert surface. Hence, by construction, genus$(K_i)\leq g$. It
follows that $g_s(K)\leq g$, where $g_s$ denotes the smooth slice
genus, and thus that $\mid s(K_i)\mid \leq 2g$  and $\mid
\tau(K_i)\mid \leq g$ by \cite{Ra} and \cite{OS} respectively.
Consider the function $ST$ from $\{K_i\}$ to
$\mathbb{Z}\times\mathbb{Z}$ given by
$ST(K_i)=(s(K_i),\tau(K_i))$. Since the image of this function is
a finite set, there is a subsequence $K_{i_j}$ on which $ST$ is
constant. Redefining our $K_i$ to be this subsequence gives the
claimed result.
\end{proof}

\section{Algebraic $n$-solution}
\label{sec:algebraic} The purpose of this section is to complete
the proof of Theorem~\ref{thm:injectivity}, which relies on
Theorem~\ref{thm:paininbutt}. In this section we define and
investigate an \emph{algebraic $n$-solution}. One might think of
this as an algebraic abstraction of an $n$-cylinder (or an
$n$-solution). In \cite{CT}, Cochran and Teichner defined an
algebraic $n$-solution. Our notion is much more general. In
particular, an algebraic $n$-solution in \cite{CT} is defined
using only the free group of rank 4, but our version is defined
using the free group of (arbitrary) even rank. Moreover, the proof
of \cite{CT} that a (geometric) $n$-solution induces an algebraic
$n$-solution is valid only for fibered knots of genus $2$. For a
group $G$, let $G_k$ denote $G/G\2k_r$ where $G\2k_r$ is the
$k$-th rational derived group of $G$. Then, as noted before, $G_k$
is a $(k-1)$-solvable PTFA group and embeds into its skew quotient
field of fractions which is denoted by $\bbk(G_k)$.

\begin{defn}
\label{defn:algebraic} Let S be a group such that $H_1(S;\bbq)\neq
0$. Suppose $F\overset{i}{\lra} S$ is a fixed homomorphism from
the free group of rank $2g$. A nontrivial homomorphism $r:S \to G$
is called an {\em algebraic $n$-solution} ($n\ge0$) (for
$F\overset{i}{\lra} S$) if the following hold :

\begin{enumerate}
\item For each $0\le k\le n-1$ the
image of the following composition,  after tensoring with the
quotient field $\bbk(G_k)$ of $\bbz G_k$, is nontrivial:
$$
H_1(S;\bbz G_k)\overset{r_*}{\lra} H_1(G;\bbz G_k)\cong
G\2k_{r}/[G\2k_{r},G\2k_{r}]\thra G\2k_{r}/G^{(k+1)}_{r}.
$$
\item For each $0\le k\le n$, the map $ H_1(F;\bbz
G_k)\overset{i_*}{\lra} H_1(S;\bbz G_k)$, after tensoring with the
quotient field $\bbk(G_k)$, is surjective.
\end{enumerate}
\end{defn}

\begin{rem}
\label{rem:algebraic}

\begin{enumerate}
\item If $n\geq 0$, then, by combining conditions $(1)$ and $(2)$
for $k=0$, we conclude that $ H_1(F;\mathbb{Q})\lra
H_1(G;\mathbb{Q})$ is non-trivial. Thus, for some $i$,
$r_*(x_{i})$ is non-trivial in $G_1=G/G^{(1)}_{r}$ where
$\{x_1,x_2,\ldots , x_m\}$ is a generating set for $F$.

\item In all the applications the image of the map in $(1)$ above
has rank at least one half the rank of $H_1(S;\bbz G_k)$.

\item If $r:S\to G$ is an algebraic $n$-solution
then, for any $k<n$ it is an algebraic $k$-solution.
\end{enumerate}
\end{rem}

\noindent Of course we need to establish that the primary
geometric examples do in fact satisfy the algebraic conditions
above.

\begin{prop}
\label{prop:geomalg} Suppose $K$ is a knot for which the degree of
the classical Alexander polynomial is greater than 2 (if $n=1$
then degree equal to $2$ is allowed). Suppose $M$ is the zero
surgery on $K$ in $S^3$, $\Sigma$ is a capped-off Seifert surface
(of genus $g$) for $K$, and $S \equiv \1(M)^{(1)}$. Suppose $F\to
\1(M-\Sigma)$  is any map inducing an isomorphism on
$H_1(F;\bbq)$. Let $i$ be the composition $F\to \1(M-\Sigma)\to
S$. Suppose W is an $n$-cylinder one of whose boundary components
is $M$. Let $G \equiv \1(W)^{(1)}$ Then the map $j:S\lra G$
(induced by inclusion) is an algebraic $n$-solution for $i:F\to
S$.
\end{prop}
\begin{proof}
First we will establish property (1) of
Definition~\ref{defn:algebraic}. Since property ($1$) is vacuous
if $n=0$, we assume $n\geq 1$. Fix an arbitrary integer $k$,
$0\leq k < n$. We need to consider the map $H_1(S;\bbz
G_k)\overset{j_*}{\lra} H_1(G;\bbz G_k)$ induced by the inclusion
map from $M$ into $W$. Let $M_\infty$ be the infinite cyclic cover
of $M$ and $W_\infty$ be the infinite cyclic cover of $W$. First
note that since $\pi_1(M_\infty)=S$ and $\pi_1(W_\infty)=G$, the
map $j_*$ is identical to $H_1(M_\infty;\mathbb{Z}
G_k)\overset{j_*}{\lra} H_1(W_\infty;\mathbb{Z} G_k)$. Now let
$\G=\1(W)/\1(W)^{(k+1)}_{r}$ so $\G$ is PTFA.  So we have the
inclusion induced map of $\mathbb{Z}\G$-modules
$j_*:H_1(M;\mathbb{Z}\G)\lra H_1(W;\mathbb{Z}\G)$. But $G_k$ is
equal to the commutator subgroup of $\G$ since $H_1(W)\cong
\mathbb{Z}$ implying that $\G_r ^{(1)}=\G ^{(1)}=G/G\2k_{r}\equiv
G_k$. Therefore this map can be also viewed as a map of
$\mathbb{Z}G_k$-modules. We claim that, as a map of $\mathbb{Z}
G_k$ modules, this is identical to our original map
$j_*:H_1(M_\infty;\mathbb{Z} G_k)\lra H_1(W_\infty;\mathbb{Z}
G_k)$. This is because $H_1(M;\mathbb{Z}\G)$ is merely the first
homology group of the total space of the $\G$-covering space of
$M$ viewed as a $\mathbb{Z}\G$-module and this total space is the
same as the total space of the $\G^{(1)}$-covering space of
$M_\infty$. Thus, as $Z\G^{(1)}$-modules, $H_1(M;\mathbb{Z}\G)$ is
the same as $H_1(M_\infty;\mathbb{Z}\G^{(1)})$. Hence we are
reduced to studying $j_*:H_1(M;\mathbb{Z}\G)\lra
H_1(W;\mathbb{Z}\G)$ as a map of $\mathbb{Z}G_k$-modules and we
need to show that it is non-trivial, even after tensoring with
$\mathbb{K}(G_k)$, i.e. after localizing with respect to the set
$R=\mathbb{Z}\G^{(1)}-\{0\}$. Since $\mathbb{Z}\G$ is an Ore
Domain, $R^{-1}\mathbb{Z}\G$ is a flat $\mathbb{Z}\G$-module and
so we just need to show that the map
$j_*:H_1(M;R^{-1}\mathbb{Z}\G)\lra H_1(W;R^{-1}\mathbb{Z}\G)$ is
nontrivial as a map of $\mathbb{K}(G_k)$ modules. For the
remainder of this proof, let $\mathbb{K}$ be the quotient field of
$\mathbb{Z}\G^{(1)}=\mathbb{Z}G_k$, which is what we have called
$\mathbb{K}(G_k)$. Finally, note that, since $H_1(\G)\cong
\mathbb{Z}$, $\mathbb{Z}\G$ can be identified with the twisted
Laurent polynomial ring $\mathbb{Z}\G^{(1)}[t^{\pm 1}]$ and,
similarly, the localized ring $R^{-1}\mathbb{Z}\G$ can be
identified with $\mathbb{K}[t^{\pm 1}]$. In summary, we are
reduced to studying $j_*:H_1(M;\mathbb{K}[t^{\pm 1}])\lra
H_1(W;\mathbb{K}[t^{\pm 1}])$ as a map of $\mathbb{K}$-modules,
and we seek to show that it is non-trivial.

We now claim that Theorem~\ref{thm:rank} applies to the map
$j_*:H_1(M;\mathbb{K}[t^{\pm 1}])\lra H_1(W;\mathbb{K}[t^{\pm
1}])$. Observe that $\G ^{n}=0$ since $\G^{k+1}=\{e\}$ and, since
$k+1 \leq n$, $\G^{n}=\{e\}$. Thus $\G$ is an $(n-1)$-solvable
PTFA group as required. We conclude that the $\mathbb{K}$-rank of
$j_*:H_1(M;\mathbb{K}[t^{\pm 1}])\lra H_1(W;\mathbb{K}[t^{\pm
1}])$ is at least $(d-2)/2$ if $n>1$ and is precisely $d/2$ if
$n=1$, where $d\equiv rk_{\mathbb{Q}}H_1(M_\infty;\mathbb{Q})$. It
is well known that $d$ is equal to the degree of the classical
Alexander polynomial of $K$ which is, by hypothesis, at least $4$
(or, if $n=1$ at least $2$). Hence in all cases the above rank is
positive and therefore we have shown that our original map ,
$H_1(S;\bbz G_k)\overset{j_*}{\lra} H_1(G;\bbz G_k)$, is
non-trivial even after tensoring with $\mathbb{K}(G_k)$.

Since the kernel of $G\2k_{r}/[G\2k_{r},G\2k_{r}]\lra
G\2k_{r}/G^{(k+1)}_{r}$ is $\mathbb{Z}$-torsion, this map is an
isomorphism after tensoring with $\mathbb{K}(G_k)$ (which contains
$\mathbb{Q})$. Combining this  with our previous conclusion, we
see that $j$ satisfies property $(1)$ of the definition of an
algebraic $n$-solution.

Now we establish property $(2)$ of an algebraic $n$-solution. We
consider an arbitrary $k$ with $0\leq k \leq n$. By flatness, what
we need to establish is the surjectivity of
$i_*:H_1(F;\mathbb{K})\lra H_1(S;\mathbb{K})$ (recall we are
sometimes abbreviating $\mathbb{K}(G_k)$ by $\mathbb{K}$). Let
$Y=M-\Sigma$ and $W$ be a wedge of $2g$ circles. By choice of $f$
there is a continuous map $W\to Y$ inducing $F\to \pi_1(Y)$ that
is $1$-connected on rational homology. It follows from
\cite[Proposition 2.10]{COT1} that this map induces a
$1$-connected map on homology with $\mathbb{K}(G_k)$ coefficients
for \emph{any} PTFA group $G_k$. Hence it is surjective. Since the
map $i_*$ factors through $H_1(Y;\mathbb{K})$, it suffices to
prove that the map $H_1(Y;\mathbb{K})\to
H_1(S;\mathbb{K})=H_1(M_\infty;\mathbb{K})$ is surjective. If $S$
were finitely generated, this would follow as above from
\cite[Proposition 2.10]{COT1}. In other words, for finitely
generated groups $S$ property $(2)$ follows from just knowing
property $(2)$ for the base case $k=0$. Unfortunately this
Proposition fails for general non-finitely generated groups and,
for us, $S$ will \emph{not} be finitely generated unless $K$ is
fibered. Also note that the rank of $H_1(M_\infty;\mathbb{K})$ may
well increase as the integer $k$ increases. Thus, to have a hope
of establishing property $(2)$, the rank of $H_1(Y;\mathbb{K})$
(which, since $Y$ has the homotopy type of a finite $2$-complex,
is bounded above, independently of $k$, by $2g$ \cite[Prop.2.10
and 2.11]{COT1}) must be a universal upper bound for the ranks of
$H_1(M_\infty;\mathbb{K})$ for \emph{any k} and \emph{any}
coefficient system. Fortunately, the technology to establish this
``universal'' upper bound essentially already exists due to work
of Shelly Harvey. We indicate how her work indeed establishes the
result we need. For intuition, first consider the case $k=0$ where
$G_k=0$ and $\mathbb{K}=\mathbb{Q}$. Then the result we claim is
that any finite generating set for the vector space
$H_1(Y;\mathbb{Q})$ generates the Alexander module of $K$ (or $M$)
as a \emph{rational vector space} (a well-known result in
classical knot theory). We proceed with the proof of the general
case. We observed above that, as $\mathbb{Z}G_k$-modules,
$H_1(M;\mathbb{Z}\G)$ is the same as $H_1(M_\infty;\mathbb{Z}G_k)$
Thus, as $\mathbb{K}$-modules, $H_1(M;\mathbb{K}[t^{\pm 1}])$ is
$H_1(M_\infty;\mathbb{K})$ (see also the discussion above
Definition~\ref{defn:localized Alexander module}). In
\cite[Proposition 7.4]{Ha} (see also \cite[section 6]{C}), Harvey
shows that $H_1(M;\mathbb{K}[t^{\pm 1}])$ has a square
presentation matrix (as a $\mathbb{K}[t^{\pm 1}]$-module) whose
entries are polynomials of degree at most $1$, with respect to
generators that are the images under inclusion of an arbitrary
generating set for $H_1(Y;\mathbb{K}[t^{\pm 1}])$. The latter is a
free module isomorphic to $H_1(Y;\mathbb{K})\otimes_{\mathbb{K}}
\mathbb{K}[t^{\pm 1}]$ (see the discussion in \cite[section
7]{Ha}). Choose a generating set $\{e_1,...,e_m\}$ for
$H_1(Y;\mathbb{K})$ as a $\mathbb{K}$-vector space. Then use the
set $\{e_i\otimes 1\}$ as a $\mathbb{K}[t^{\pm 1}]$ generating set
for the module $H_1(Y;\mathbb{K}[t^{\pm 1}])$. Harvey's
presentation result shows that these elements generate
$H_1(M;\mathbb{K}[t^{\pm 1}])$ \emph{as a $\mathbb{K}[t^{\pm
1}]$-module}. We claim that a closer analysis of some other work
of Harvey shows the stronger fact that these actually generate as
a $\mathbb{K}$-module! This will finish the verification that the
map $H_1(Y;\mathbb{K})\to H_1(M_\infty;\mathbb{K})$ is surjective
and thus finish the verification of property $(2)$. In
\cite[Proposition 9.1]{Ha}, Harvey shows that any such matrix as
above (whose entries are polynomials of degree at most $1$)
presents a module of rank $m$ over $\mathbb{K}$. It is only
necessary to examine her proof carefully to see that this stronger
statement is true: The given set of $m$ generators is a
$\mathbb{K}$-generating set for the module. Harvey's proof
involves changing the presentation matrix by certain allowable
matrix operations and finally arriving at a simple matrix which
she explicitly shows satisfies this stronger statement. Therefore
it is only necessary to check that all the matrix operations she
uses do in fact preserve the veracity of this statement. For
example, clearly $\{x_1,...,x_m\}$ is a $\mathbb{K}$-generating
set for the (right) module if and only if
$\{x_1-(x_2)k,x_2...,x_m\}$ is also, where $k$ is any non-zero
element of $\mathbb{K}$. This translates into the fact that adding
a non-zero left $\mathbb{K}$-multiple of a row of Harvey's matrix
to another row is an allowable operation for our purposes. The key
point is to \emph{not} allow a $\mathbb{K}[t^{\pm 1}]$-multiple.
In fact the only matrix operation from Harvey's list that could
lead to problems is : add to any row a left $\mathbb{K}[t^{\pm
1}]$-linear combination of the other rows (because this
corresponds to a change of generators). One only needs to notice
that, in fact, in her proof she never uses the full generality of
this operation. She only adds to any row a left multiple of
another row by a non-zero element of $\mathbb{K}$ (not
$\mathbb{K}[t^{\pm 1}]$). This completes the proof.
\end{proof}

The following theorem greatly generalizes \cite[Theorem 6.3]{CT}
because our definition of an algebraic $n$-solution is much more
general. Our proof follows the proof of \cite[Theorem 6.3]{CT}
closely. In doing so we had to make decisions about what to
include and what to reference. Because the result is already
(logically) extremely demanding on the reader, for the reader's
convenience, we have erred on the side of duplicating material and
have included a complete proof. In addition, our proof is simpler
in certain places.

\begin{thm}
\label{thm:paininbutt} Suppose we are given $F\overset{i}{\lra} S$
as in Definition~\ref{defn:algebraic}. For each $n\ge0$ there is a
finite collection $\SP_n$ of sets consisting of $2g-1$, if $n>0$,
or $2g$, if $n=0$, elements of $F\2n$ (we refer to such a set as a
tuple even though it is unordered), with the following property:
For any algebraic $n$-solution $r$ for $F\overset{i}{\lra} S$, at
least one such tuple (which will be called a {\bf special tuple}
for $r$) maps to a generating set under the composition:
$$
F\2n\lra S\2n/S\nn\cong H_1(S;\bbz S_n)\overset{r_*}{\lra}
H_1(S;\bbz G_n)\lra H_1(S;\bbz G_n)\otimes_{\bbz G_n} \bbk(G_n)
$$
where $\bbk(G_n)$ is the skew quotient field of fractions of $\bbz
G_n$.
\end{thm}

\begin{proof} We remark that to be used for the proof of the main
theorem, it is very important to define the collections $\SP_n$ so
as to depend only on the knot $K$. In particular, they must not
depend on the existence of any particular $n$-cylinder.

Let $\{x_1,...,x_{2g}\}$ be a generating set of the free group
$F$. Set $\SP_0 \equiv \{\{x_1,...,x_{2g}\}\}$, the collection
consisting of a single $2g$-tuple. Set $\SP_1 \equiv
\{\{[x_i,x_1],\dots,
[x_i,x_{i-1}],[x_i,x_{i+1}],\ldots,[x_i,x_{2g}]\}\,\,|\,\,1\le i
\le 2g\}$, the collection consisting of $(2g-1)\cdot (2g)$ number
of $(2g-1)$-tuples. Supposing $\SP_k$ $(k\ge 1)$ has been defined,
define $\SP_{k+1}$ recursively as follows. For each
$\{w_1,...,w_{2g-1}\}\in$ $\SP_k$ include the $2g-1$-tuple
$\{z_1,...,z_{2g-1}\}$ in $\SP_{k+1}$ if $z_i=[w_i,w_i^{x_j}]$
($1\le i \le 2g-1, 1\le j \le 2g$, and $i\ne j$)  or if
$z_i=[w_i,w_k]$ for some $1\le i,k\leq 2g-1$ and $i\ne k$. Here
$w_i^{x_j} \equiv x_j^{-1}w_i x_j$. Clearly $z_i\in F^{(k+1)}$.

Now we fix $n$ and show $\SP_n$ satisfies the conditions of the
theorem. Fix an algebraic $n$-solution $r:S\lra G$. We must show
that there exists a special tuple in $\SP_n$ corresponding to $r$.
This is trivially true for $n=0$ using property $(2)$ of the
definition of an algebraic $n$-solution, so we assume $n\ge1$. Now
we need some preliminary definitions.

Recall that $F$ is the free group on $\{x_1,\dots,x_{2g}\}$. Its
classifying space has a standard cell structure as a wedge of $2g$
circles $W$. Our convention is to consider its universal cover
$\wt W$ as a right $F$-space as follows. Choose a preimage of the
$0$-cell as basepoint denoted $*$. For each element $w\in
F\equiv\1(W)$, lift $w^{-1}$ to a path $(\tl w^{-1})$ beginning at
$*$. There is a unique deck translation $\Phi(w)$ of $\wt W$ which
sends $*$ to the endpoint of this lift. Then $w$ acts on $\wt W$
by $\Phi(w)$. This is the conjugate action of the usual left
action as in \cite{Ma}. Taking the induced cell structure on $\wt
W$ and tensoring with an arbitrary left $\bbz F$-module $A$ gives
an exact sequence
\begin{equation}
0\lra H_1(F;A)\overset{d}{\lra} A^{2g}\lra A\lra H_0(F;A)\lra0.
\end{equation}
Specifically consider $A=\bbz G$ where $\bbz F$ acts by left
multiplication via a homomorphism $\phi:F\lra G$. From the
interpretation of $H_1(F;\bbz G)$ as $H_1$ of a $G$-cover of $W$,
one sees that an element $g$ of $\Ker(\phi)$ can be considered as
an element of $H_1(F;\bbz G)$. We claim that the composition
$\Ker(\phi)\lra H_1(F;\bbz G)\overset{d}{\lra}(\bbz G)^{2g}$ can
be calculated using the ``free differential calculus''
$\p=(\p_1,\dots,\p_{2g})$ where $\p_i:F\lra\bbz F$. Specifically
we assert that the diagram below commutes
$$
\begin{diagram}\dgARROWLENGTH=1.2em
\node{F} \arrow[2]{e,t}{\partial} \node[2]{(\bbz F)^{2g}}
\arrow{s,r}{\phi}\\
\node{\Ker(\phi)} \arrow{n} \arrow{e} \node{H_1(F;\bbz G)}
\arrow{e,t}{d} \node{(\bbz G)^{2g}}
\end{diagram}
$$
where $\p_i(x_j)=\d_{ij}$, $\p_i(e)=0$ and
$\p_i(gh)=\p_ig+(\p_ih)g^{-1}$ for each $1\le i\le{2g}$. Note that the
usual formula for the standard left action is
$d_i(gh)=d_i(g)+gd_i(h)$. Our formula is obtained by setting
$\p_i=\bar d_i$ where $^-$ is the involution on the group ring.
This is justified in more detail in \cite[Section 6]{CT}.

Henceforth we abbreviate maps of the form $(r,...,r):(\bbz
F_n)^{2g}\lra(\bbz G_n)^{2g}$ as $r$. Note that $r\circ\pi_k:F\to
F_k\to G_k$ is the same as $\pi_k\circ r:F\to S\to G\to G_k$.

\begin{lem}
\label{lem:good} Given an algebraic $n$-solution $r:S\lra G$ for
$F\to S$, for each $k$, $1\le k\le n$ there is an ordering of the
basis elements $\{x_1,...,x_{2g}\}$ of $F$ such that there exists
at least one tuple $\{w_1,...,w_{2g-1}\}$$\in\SP_k$ with the
following {\bf good} property: The set of 2g-1 vectors (obtained
as $i$ varies from $1$ to $2g-1$)
$(r\pi_k\p_1w_i,...,r\pi_k\p_{2g-1}w_i)$ (i.e. the vectors
consisting of the first $2g-1$ coordinates of the images of the
$w_i$ under the composition $F\2k\overset{\pi_k\p}{\lra}(\bbz
F_k)^{2g}\overset{r}{\lra}(\bbz G_k)^{2g})$ is right linearly
independent over $\bbz G_k$.
\end{lem}

\begin{proof}[Proof that Lemma~\ref{lem:good} $\Longrightarrow$
Theorem~\ref{thm:paininbutt}] The set $\SP_n$ has been defined.
Recall that we are assuming that $n\geq 1$. Given an algebraic
$n$-solution the Lemma provides a tuple
$\{w_1,...,w_{2g-1}\}$$\in\SP_n$ which has the {\bf good}
property. We verify that any \textbf{good} tuple
$\{w_1,...,w_{2g-1}\}$ for $r$ is a {\bf special} tuple for $r$.
Consider the diagram below. Recall $G_n=G/G\2n_{r}$.
$$
\begin{diagram}\dgARROWLENGTH=1.2em
\node{F\2n} \arrow{e} \node{H_1(F;\bbz F_n)} \arrow{e,t}{d}
\arrow{s,r}{r_*} \node{(\bbz F_n)^{2g}} \arrow{s,r}{r}\\
\node[2]{H_1(F;\bbz G_n)} \arrow{e,t}{d'} \arrow{s,r}{i_*}
\node{(\bbz G_n)^{2g}}\\
\node[2]{H_1(S;\bbz G_n)}
\end{diagram}
$$

The horizontal composition on top is $\pi_n\circ\p$. The right-top
square commutes by naturality of the sequence (1) above. The rank
of $H_1(F;\bbz G_n)$ over $\bbk(G_n)$ is $2g-1$ by \cite[Lemma
3.9]{C}. Since $d'$ is a monomorphism, and $i_*$ is an epimorphism
after tensoring with $\bbk(G_n)$, to show $\{w_1,...,w_{2g-1}\}$
is special, it suffices to show that the set
$\{r\pi_n\p(w_1),...,r\pi_n\p(w_{2g-1})\}$ is $\bbz G_n$-linearly
independent in $(\bbz G_n)^{2g}$. This follows immediately from
the {\bf good} property of the tuple $\{w_1,...,w_{2g-1}\}$.
This completes the verification that the Lemma implies the
Theorem.
\end{proof}

\begin{proof}[Proof of Lemma~\ref{lem:good}] The integer $n$
is fixed throughout. Let $r:S\to G$ be a fixed algebraic
$n$-solution. Recall that, by Remark~\ref{rem:algebraic}(1), for
any ordering of the basis elements of $F$, there is some $i$ such
that $r_*(x_{i})$ is non-trivial in $G_1=G/G^{(1)}_{r}$. Thus by
reordering we may assume that $r_*(x_{2g})$ is non-trivial in
$H_1(G;\bbq)=G_1\ox\bbq$. Since $G_1\equiv G/G\21_{r}$ is
torsion-free, this means that $r\1(x_{2g})$ is non-trivial in
$G_1$.

We will prove the Lemma by induction on $k$. We begin with $k=1$.
Consider the $2g-1$-tuple
$\{z_1,...,z_{2g-1}\}$$=\{[x_{2g},x_1],...,[x_{2g},x_{2g-1}]\}\in\SP_1$.
We claim that it has the {\bf good} property.
Using $d$ as shorthand for $\p_j$, one computes that
$d([g,h])=dg+(dh)g^{-1}-(dg)gh^{-1}g^{-1}-(dh)hgh^{-1}g^{-1}$
Using this we compute that, for any $j$, $1\leq j\leq 2g-1$,
$\p_jz_i$ is $x_{2g}^{-1}-[x_i,x_{2g}]$ if $j=i$ and is otherwise
zero. Thus, only the $i^{th}$ coordinate of
$(r\pi_k\p_1z_i,...,r\pi_k\p_{2g-1}z_i)$ is (possibly) non-zero.
Therefore the square matrix whose columns are these vectors is a
diagonal matrix and to establish the {\bf good} property, it
suffices to show that $r\1(x_{2g}^{-1}-[x_i,x_{2g}])\neq 0$ since
this certainly implies that $r\pi_k(x_{2g}^{-1}-[x_i,x_{2g}])\neq
0$. This follows since $r\1([x_i,x_{2g}])=e$ and
$r\1(x_{2g}^{-1})$ is non-trivial by assumption. Therefore the
base of the induction $(k=1)$ is established.

Now suppose the conclusions of the Lemma have been established for
$1,\dots,k$ where $k<n$. We establish them for $k+1$. Note that
$r:S\to G$ is also an algebraic $k$-solution for $i:F\to S$. By
induction the Lemma holds for $k<n$, so there is a tuple
$\{w_1,...,w_{2g-1}\}$$\in\SP_k$ that has the {\bf good} property.
First, we claim that there is at least one of the $w_i$ (which by
relabelling we will assume is $w_1$)
 such that $r\pi_{k+1}(w_1)\neq e$ in
$G_{k+1}$. For, by the \emph{proof} of Lemma~\ref{lem:good}
$\Longrightarrow$ Theorem~\ref{thm:paininbutt}, we know that
$\{w_1,...,w_{2g-1}\}$$\in\SP_k$ is \textbf{special} for $r$. Thus
under the composition $ F\2k\lra H_1(S;\bbz
S_k)\overset{r_*}{\lra}H_1(S;\bbz G_k){\lra}H_1(S;\bbz G_k)\otimes
\bbk(G_k)$ the set $\{w_i\}$ maps to a generating set. Combined
with the fact that $r$ is also an algebraic $n$-solution, we see
that the composition of the above with the map
$$
H_1(S;\bbz G_k)\otimes \bbk(G_k) \overset{r_*}{\lra}H_1(G;\bbz
G_k)\otimes \bbk(G_k) \cong G\2k_{r}/G^{(k+1)}_{r}\otimes
\bbk(G_k)
$$
is nontrivial when restricted to $\{w_i\}$. On the other hand the
combined map $F\2k\lra G\2k_{r}/G^{(k+1)}_{r}$ is clearly given by
$w_i\mapsto r\pi_{k+1}(w_i)$ so it is not possible that all of the
elements $r\pi_{k+1}(w_i)$ lie in $G^{(k+1)}_{r}$. Hence we may
assume that $r\pi_{k+1}(w_1)\neq e$ in $G_{k+1}$.

Consider the tuple $\{z_1,...,z_{2g-1}\}\in \SP_{k+1}$,  where
$z_i=[w_i,w_i^{x_{2g}}]$ if $r\pi_{k+1}(w_i)\neq e$ in $G_{k+1}$,
and $z_i=[w_i,w_1]$ if $r\pi_{k+1}(w_i)= e$ in $G_{k+1}$. We will
show that this tuple has the {\bf good} property, finishing the
inductive proof of Lemma~\ref{lem:good}.

For the remainder of this proof we write $x$ for $x_{2g}$,
suppressing the subscript. We need to show that the set of $2g-1$
vectors (obtained as $i$ varies from $1$ to $2g-1$)
$(r\pi_{k+1}\p_1z_i,...,r\pi_{k+1}\p_{2g-1}z_i)$ is $\bbz
G_{k+1}$-linearly independent. For this purpose we compute, for
each $i$, $\p_jz_i$. This computation falls into two cases.

\noindent{\bf Case 1}: $r\pi_{k+1}(w_i) \neq e$.

Let $d$ be shorthand for  $\p_j$ for $1\leq j \leq 2g-1$. One has
$dx=0$, $d(g^{-1})=(dg)g$, and
$dg^x=(dg)x$. Using these one computes that
$d([w_i,w_i^x])=(dw_i)p_i$ where $p_i$ is independent of $j$ and
is equal to $1+xw_i^{-1}-(w_i^x)^{-1}[w_i^x,w_i]-x[w_i^x,w_i]$.
Thus for any value of $i$ that falls under Case $1$, the vector
$(\p_1z_i,...\p_{2g-1}z_i)$ is  a right multiple of the vector
$(\p_1w_i,...\p_{2g-1}w_i)$ by the element $p_i$. Hence the vector
in question, $(r\pi_{k+1}\p_1z_i,...,r\pi_{k+1}\p_{2g-1}z_i)$ is a
right multiple of the vector
$(r\pi_{k+1}\p_1w_i,...,r\pi_{k+1}\p_{2g-1}w_i)$ by the element
$r\pi_{k+1}p_i$ which lies in $\bbz G_{k+1}$. The right factor
$r\pi_{k+1}p_i$ is seen to be non-trivial in $\bbz G_{k+1}$ as
follows. Note $r\pi_{k+1}p_i$ is a linear combination of $4$ group
elements $e$, $r\pi_{k+1}(xw_i^{-1})$, $r\pi_{k+1}((w_i^x)^{-1})$,
and $r\pi_{k+1}(x)$ in $G_{k+1}$. For $r\pi_{k+1}p_i$ to vanish in
$\bbz G_{k+1}$ the group elements would have to pair up in a
precise way and in particular in such a way that
$r\pi_{k+1}(x)=r\pi_{k+1}(xw_i^{-1})$ in $G_{k+1}$. This is a
contradiction since $r\pi_{k+1}(w_i)\neq e$ by hypothesis. No
other pairing is possible because the projections of the $4$
elements to $G_1$ are $e$, $r\1(x)$, $e$ and $r\1(x)$ and we have
already noted that $r\1(x)=r\1(x_{2g})$ is non-trivial in $G_1$.
Since these right factors are non-trivial and since $\bbz G_{k+1}$
has no zero divisors, the right linear independence of the
collection of vectors obtained by ignoring the $p_i$,
$(r\pi_{k+1}\p_1w_i,...,r\pi_{k+1}\p_{2g-1}w_i)$,
 would be sufficient to imply the right linear independence of the
original set of vectors. We denote these new vectors by
$\mathbf{v_i^{k+1}}$. So far we have only dealt with those values
of $i$ that fall under Case $1$.

\noindent{\bf Case 2}: $r\pi_{k+1}(w_i)=e$.

 As above one computes that $d([w_i,w_1])=(dw_i)q_i +
 (dw_1)(w_i^{-1}-[w_1,w_i])$
where $q_i$ is independent of $j$ and is equal to
$1-w_1^{-1}[w_1,w_i]$. Note that under the map $r\pi_{k+1}$ the
factor $(w_i^{-1}-[w_1,w_i])$ goes to zero. Thus for any value of
$i$ that falls under Case $2$, the vector in question,
$(r\pi_{k+1}\p_1z_i,...,r\pi_{k+1}\p_{2g-1}z_i)$ is a right
multiple of the vector
$(r\pi_{k+1}\p_1w_i,...,r\pi_{k+1}\p_{2g-1}w_i)$ by the element
$r\pi_{k+1}q_i$. We denote the latter vector by
$\mathbf{v_i^{k+1}}$ as above. An argument just as in Case $1$
shows that the right factor $r \pi_{k+1}q_i$ is non-trivial, using
the nontriviality of $r\pi_{k+1}(w_1)$.

Now our objective is to show that the set of vectors
$\mathbf{v_i^{k+1}}$ is $\bbz G_{k+1}$-linearly independent.
Recall that our hypothesis is that the set of vectors
$(r\pi_{k}\p_1w_i,...,r\pi_{k}\p_{2g-1}w_i)$, which we denote by
$\mathbf{v_i^{k}}$, is $\bbz G_{k}$-linearly independent. Note
that, for each $i$, $\mathbf{v_i^{k}}$ is the image of
$\mathbf{v_i^{k+1}}$ under the canonical projection $(\bbz
G_{k+1})^{2g-1}\lra (\bbz G_k)^{2g-1}$. We assert that the linear
independence of $\mathbf{v_i^{k}}$ implies the linear independence
of $\mathbf{v_i^{k+1}}$ since the kernel of $G_{k+1}\lra G_k$ is a
torsion free abelian group and hence a $D(\mathbf{Z})$-group in
the sense of R. Strebel. Details follow. This will complete the
verification, that $\{z_1,...z_{2g-1}\}$ has the {\bf good}
property.

To establish our assertion above, consider that the vectors
$\mathbf{v_i^{k+1}}$ describe an endomorphism $f$ of free right
$\bbz G_{k+1}$ modules $f:(\bbz G_{k+1})^{2g-1}\lra(\bbz
G_{k+1})^{2g-1}$ given by multiplying on the left by the matrix
$M$ whose columns are the $\mathbf{v_i^{k+1}}$. Then $f$ induces
such a map $\bar f:(\bbz G_k)^{2g-1}\lra(\bbz G_k)^{2g-1}$ given
by the matrix $\bar M$ obtained by projecting each entry of $M$.
Thus the columns are the $\mathbf{v_i^{k}}$. Let
$H=G\2k_{r}/G^{(k+1)}_{r}$ and note that $\bbz
G_{k+1}\ox_{\bbz[H]}\bbz\cong\bbz G_k$ where $a\ox1\mapsto\bar a$,
as $\bbz G_{k+1}-\bbz$ bimodules. Moreover (under this
identification) $f$ descends to $f\ox\id:(\bbz
G_k)^{2g-1}\lra(\bbz G_k)^{2g-1}$ sending $\bar v$ (for $v\in(\bbz
G_{k+1})^{2g-1})$ to $\bar f (\bar v)$, thus agreeing with $\bar
f$ above. Our hypothesis on $\mathbf{v_i^{k}}$ guarantees that the
columns of $\bar M$ are right linearly independent and hence that
$\bar f$ is injective. Since $H$ is torsion-free-abelian, a
theorem of Strebel \cite[Section 1]{Str} ensures that the
injectivity of $f\ox\id=\bar f$ implies the injectivity of $f$.
This in turn implies the linear independence of the columns of
$M$, which are the $\mathbf{v_i^{k+1}}$.

\end{proof}
\end{proof}


\begin{thebibliography}{COT1}
\bibitem[CS]{CS} S.~Cappell and J.~L.~Shaneson, {\em The codimension
two placement problem and homology equivalent manifolds}, Ann.~of
Math.~(2) \textbf{99} (1974), 277--348.

\bibitem[CG]{CG} A.~J.~Casson and C.~McA.~Gordon, {\em Cobordism of
classical knots}, printed notes, Orsay, 1975, published in A la
recherche de la topologie perdue 62, ed. Guilllou and Marin,
Progress in Mathematics, 1986.

\bibitem[ChG]{ChG}
J.~Cheeger, M.~Gromov, {\em Bounds on the von Neumann dimension of
$L\sp 2$-cohomology and the Gauss-Bonnet theorem for open
manifolds}, J.~Differential Geom.~\textbf{21} (1985), no.~1,
1--34.


\bibitem[COT1]{COT1} T.~D.~Cochran, K.~E.~Orr, and P.~Teichner,
{\em Knot concodance, Whitney towers and $L^2-$signatures},
Ann.~of Math. (2) \textbf{157} (2003), no.~2, 433--519.

\bibitem[COT2]{COT2} T.~D.~Cochran, K.~E.~Orr, and P.~Teichner,
{\em Structure in the classical knot concordance group},
Comment.~Math.~Helv.~\textbf{79} (2004), no. 1, 105--123.

\bibitem[C]{C} T.~Cochran, {\em Noncommutative Knot Theory},
Algebr.~Geom.~Topol. \textbf{4} (2004), 347--398.


\bibitem[CT]{CT} T.~D.~Cochran, P.~Teichner, {\em Knot concordance and von Neumann $\eta$-invariants}, preprint, 2004,
arXiv:math.~GT/0411057.


\bibitem[F]{F}
M.~H.~Freedman, {\em The topology of four dimensional manifolds},
J.~Differential Geom. \textbf{17} (1982), no.~3, 357--453.

\bibitem[FQ]{FQ}
M.~H.~Freedman and F.~Quinn, \emph{Topology of 4-manifolds},
Princeton Mathematical Series, \textbf{39}, Princeton University
Press, Princeton, NJ, 1990.


\bibitem[Fr]{Fr} S.~Friedl, \emph{Eta invariants
as sliceness obstructions and their relation to Casson-Gordon
invariants}, Algebr.~Geom.~Topol. \textbf{4} (2004), 893--934.

\bibitem[FT]{FT} S.~Friedl and P.~Teichner, \emph{New
examples of topologically slice knots}, in preparation.

\bibitem[G]{G} P.~M.~Gilmer, {\em Slice knots in $S^3$},
Quart.~J.~Math. Oxford Ser. (2) \textbf{34} (1983), no.~135,
305--322.


\bibitem[Ha]{Ha} S.~Harvey, {\em Higher-Order Polynomial
Invariants of 3-manifolds giving lower bounds for the Thurston
Norm}, to appear in Topology, arXiv:math.GT/0307107.

\bibitem[K]{K} T.~Kim, {\em Infinite family of non-concordant
knots having the same Seifert form}, to appear in
Comment.~Math.~Helv., arXiv:math.GT/0402425.

\bibitem[KL]{KL} P.~Kirk and C.~Livingston, \emph{Twisted Alexander
invariants, Reidemeister torsion, and Casson-Gordon invariants},
Topology \textbf{38} (1999), no.~3, 635--661.

\bibitem[L]{L}
J.~Levine, \emph{Knot cobordism groups in codimension two},
Comment.~Math.~Helv. \textbf{44} (1969), 229--244.

\bibitem[Le]{Le} C.~F.~Letsche, {\em An obstruction to slicing
knots using the eta invariant}, Math.~Proc.~Cambridge Philos.~Soc.
\textbf{128} (2000), no.~2, 301--319.

\bibitem[Li]{Li}
C.~Livingston, \emph{Seifert forms and concordance}, Geom.~Topol.
\textbf{6} (2002), 403--408.

\bibitem[Ma]{Ma} W.S. Massey, {\it Algebraic Topology : An
Introduction}, Harcourt, Brace, and World, 1967.

\bibitem[OS]{OS} P.~Ozsv{\'a}th and Z.~Szab{\'o}, {\em Knot Floer
homology and the four-ball genus}, Geom.~Topol. \textbf{7} (2003),
615--639 (electronic).

\bibitem[R]{R} M.~Ramachandran, \emph{von Neumann index
theorems for manifolds with boundary}, J.~Differential Geom.
\textbf{38} (1993), no. 2, 315--349.

\bibitem[Ra]{Ra} J.~A.~Rasmussen, \emph{Khovanov homology
and the slice genus}, preprint, 2004, arXiv:math.~GT/0402131.

\bibitem[Ste]{Ste} B.~Stenstr\"om, {\em Rings of Quotients},
Springer-Verlag, 1975, New York.

\bibitem[Str]{Str} R.~Strebel, {\em Homological methods
applied to the derived series of groups}, Comment.~Math.~Helv.
\textbf{49} (1974), 302--332.

\bibitem[Wa]{Wa} C.~T.~C.~Wall, {\em Surgery on Compact Manifolds},
London Math.~Soc.~Monographs~1, Academic Press~1970.

\end{thebibliography}
\end{document}